
\input amstex



\def\spaces{\space\space\space\space\space\space\space\space\space\space}
\def\spacess{\message{\spaces\spaces\spaces\spaces\spaces\spaces\spaces}}
\spacess
\spacess
\message{Annals of Mathematics Style: Current Version: 1.1. June 10, 1992}
\spacess
\spacess

\catcode`\@=11

\hyphenation{acad-e-my acad-e-mies af-ter-thought anom-aly anom-alies
an-ti-deriv-a-tive an-tin-o-my an-tin-o-mies apoth-e-o-ses apoth-e-o-sis
ap-pen-dix ar-che-typ-al as-sign-a-ble as-sist-ant-ship as-ymp-tot-ic
asyn-chro-nous at-trib-uted at-trib-ut-able bank-rupt bank-rupt-cy
bi-dif-fer-en-tial blue-print busier busiest cat-a-stroph-ic
cat-a-stroph-i-cally con-gress cross-hatched data-base de-fin-i-tive
de-riv-a-tive dis-trib-ute dri-ver dri-vers eco-nom-ics econ-o-mist
elit-ist equi-vari-ant ex-quis-ite ex-tra-or-di-nary flow-chart
for-mi-da-ble forth-right friv-o-lous ge-o-des-ic ge-o-det-ic geo-met-ric
griev-ance griev-ous griev-ous-ly hexa-dec-i-mal ho-lo-no-my ho-mo-thetic
ideals idio-syn-crasy in-fin-ite-ly in-fin-i-tes-i-mal ir-rev-o-ca-ble
key-stroke lam-en-ta-ble light-weight mal-a-prop-ism man-u-script
mar-gin-al meta-bol-ic me-tab-o-lism meta-lan-guage me-trop-o-lis
met-ro-pol-i-tan mi-nut-est mol-e-cule mono-chrome mono-pole mo-nop-oly
mono-spline mo-not-o-nous mul-ti-fac-eted mul-ti-plic-able non-euclid-ean
non-iso-mor-phic non-smooth par-a-digm par-a-bol-ic pa-rab-o-loid
pa-ram-e-trize para-mount pen-ta-gon phe-nom-e-non post-script pre-am-ble
pro-ce-dur-al pro-hib-i-tive pro-hib-i-tive-ly pseu-do-dif-fer-en-tial
pseu-do-fi-nite pseu-do-nym qua-drat-ics quad-ra-ture qua-si-smooth
qua-si-sta-tion-ary qua-si-tri-an-gu-lar quin-tes-sence quin-tes-sen-tial
re-arrange-ment rec-tan-gle ret-ri-bu-tion retro-fit retro-fit-ted
right-eous right-eous-ness ro-bot ro-bot-ics sched-ul-ing se-mes-ter
semi-def-i-nite semi-ho-mo-thet-ic set-up se-vere-ly side-step sov-er-eign
spe-cious spher-oid spher-oid-al star-tling star-tling-ly
sta-tis-tics sto-chas-tic straight-est strange-ness strat-a-gem strong-hold
sum-ma-ble symp-to-matic syn-chro-nous topo-graph-i-cal tra-vers-a-ble
tra-ver-sal tra-ver-sals treach-ery turn-around un-at-tached un-err-ing-ly
white-space wide-spread wing-spread wretch-ed wretch-ed-ly Brown-ian
Eng-lish Euler-ian Feb-ru-ary Gauss-ian Grothen-dieck Hamil-ton-ian
Her-mit-ian Jan-u-ary Japan-ese Kor-te-weg Le-gendre Lip-schitz
Lip-schitz-ian Mar-kov-ian Noe-ther-ian No-vem-ber Rie-mann-ian
Schwarz-schild Sep-tem-ber Za-mo-lod-chi-kov Knizh-nik quan-tum Op-dam
Mac-do-nald Ca-lo-ge-ro Su-ther-land Mo-ser Ol-sha-net-sky  Pe-re-lo-mov
in-de-pen-dent ope-ra-tors
}

\Invalid@\nofrills
\Invalid@\usualspace
\newif\ifnofrills@
\def\nofrills@#1#2{\relaxnext@
  \DN@{\ifx\next\nofrills
    \nofrills@true\let#2\relax\DN@\nofrills{\nextii@}%
  \else
    \nofrills@false\def#2{#1}\let\next@\nextii@\fi
\next@}}
\def\usualspace@#1{\ifnofrills@\def\usualspace{#1}\fi}
\def\addto#1#2{\csname \expandafter\eat@\string#1@\endcsname
  \expandafter{\the\csname \expandafter\eat@\string#1@\endcsname#2}}
\newdimen\bigsize@
\def\big@#1#2{{\hbox{$\left#2\vcenter to#1\bigsize@{}%
  \right.\nulldelimiterspace\z@\m@th$}}}
\def\big{\big@\@ne}
\def\Big{\big@{1.5}}
\def\bigg{\big@\tw@}
\def\Bigg{\big@{2.5}}
\def\raggedcenter@{\leftskip\z@ plus.4\hsize \rightskip\leftskip
 \parfillskip\z@ \parindent\z@ \spaceskip.3333em \xspaceskip.5em
 \pretolerance9999\tolerance9999 \exhyphenpenalty\@M
 \hyphenpenalty\@M \let\\\linebreak}
\def\upperspecialchars{\def\ss{SS}\let\i=I\let\j=J\let\ae\AE\let\oe\OE
  \let\o\O\let\aa\AA\let\l\L}
\def\uppercasetext@#1{%
  {\spaceskip1.2\fontdimen2\the\font plus1.2\fontdimen3\the\font
   \upperspecialchars\uctext@#1$\m@th\aftergroup\eat@$}}
\def\uctext@#1$#2${\endash@#1-\endash@$#2$\uctext@}
\def\endash@#1-#2\endash@{%
\uppercase{#1}\if\notempty{#2}--\endash@#2\endash@\fi}
\def\runaway@#1{\DN@{#1}\ifx\envir@\next@
  \Err@{You seem to have a missing or misspelled \string\end#1 ...}%
  \let\envir@\empty\fi}
\newif\iftemp@
\def\notempty#1{TT\fi\def\test@{#1}\ifx\test@\empty\temp@false
  \else\temp@true\fi \iftemp@}

\font@\tensmc=cmcsc10
\font@\sevenex=cmex7
\font@\sevenit=cmti7
\font@\eightrm=cmr8 
\font@\sixrm=cmr6 
\font@\eighti=cmmi8     \skewchar\eighti='177 
\font@\sixi=cmmi6       \skewchar\sixi='177   
\font@\eightsy=cmsy8    \skewchar\eightsy='60 
\font@\sixsy=cmsy6      \skewchar\sixsy='60   
\font@\eightex=cmex8 %
\font@\eightbf=cmbx8 
\font@\sixbf=cmbx6   
\font@\eightit=cmti8 
\font@\eightsl=cmsl8 
\font@\eightsmc=cmcsc10
\font@\eighttt=cmtt8 

\loadmsam
\loadmsbm
\loadeufm
\UseAMSsymbols

\def\penaltyandskip@#1#2{\relax\ifdim\lastskip<#2\relax\removelastskip
      \ifnum#1=\z@\else\penalty@#1\relax\fi\vskip#2%
  \else\ifnum#1=\z@\else\penalty@#1\relax\fi\fi}
\def\nobreak{\penalty\@M
  \ifvmode\def\penalty@{\let\penalty@\penalty\count@@@}%
  \everypar{\let\penalty@\penalty\everypar{}}\fi}
\let\penalty@\penalty

\def\block{\RIfMIfI@\nondmatherr@\block\fi
       \else\ifvmode\vskip\abovedisplayskip\noindent\fi
        $$\def\endblock{\par\egroup$$}\fi
  \vbox\bgroup\advance\hsize-2\indenti\noindent}
\def\endblock{\par\egroup}

\def\logo@{\baselineskip2pc \hbox to\hsize{\hfil\eightpoint Typeset by
 \AmSTeX}}




\font\elevensc=cmcsc10 scaled\magstephalf
\font\tensc=cmcsc10

\font\eightsc=cmcsc10 scaled800

\font\elevenrm=cmr10 scaled \magstephalf
\font\ninerm=cmr9
\font\eightrm=cmr8
\font\sixrm=cmr6
\font\fiverm=cmr5

\font\eleveni=cmmi10 scaled\magstephalf
\font\ninei=cmmi9
\font\eighti=cmmi8
\font\sixi=cmmi6
\font\fivei=cmmi5
\skewchar\ninei='177 \skewchar\eighti='177 \skewchar\sixi='177
\skewchar\eleveni='177

\font\elevensy=cmsy10 scaled\magstephalf
\font\ninesy=cmsy9
\font\eightsy=cmsy8
\font\sixsy=cmsy6
\font\fivesy=cmsy5
\skewchar\ninesy='60 \skewchar\eightsy='60 \skewchar\sixsy='60
\skewchar\elevensy'60

\font\eighteenbf=cmbx10 scaled\magstep3

\font\twelvebf=cmbx10 scaled \magstep1
\font\elevenbf=cmbx10 scaled \magstephalf
\font\tenbf=cmbx10
\font\ninebf=cmbx9
\font\eightbf=cmbx8
\font\sixbf=cmbx6
\font\fivebf=cmbx5

\font\elevenit=cmti10 scaled\magstephalf
\font\nineit=cmti9
\font\eightit=cmti8

\font\eighteenmib=cmmib10 scaled \magstep3
\font\twelvemib=cmmib10 scaled \magstep1
\font\elevenmib=cmmib10 scaled\magstephalf
\font\tenmib=cmmib10
\font\eightmib=cmmib10 scaled 800
\font\sixmib=cmmib10 scaled 600

\font\eighteensyb=cmbsy10 scaled \magstep3
\font\twelvesyb=cmbsy10 scaled \magstep1
\font\elevensyb=cmbsy10 scaled \magstephalf
\font\tensyb=cmbsy10
\font\eightsyb=cmbsy10 scaled 800
\font\sixsyb=cmbsy10 scaled 600

\font\elevenex=cmex10 scaled \magstephalf
\font\tenex=cmex10
\font\eighteenex=cmex10 scaled \magstep3


\def\elevenpoint{\def\rm{\fam0\elevenrm}%
  \textfont0=\elevenrm \scriptfont0=\eightrm \scriptscriptfont0=\sixrm
  \textfont1=\eleveni \scriptfont1=\eighti \scriptscriptfont1=\sixi
  \textfont2=\elevensy \scriptfont2=\eightsy \scriptscriptfont2=\sixsy
  \textfont3=\elevenex \scriptfont3=\tenex \scriptscriptfont3=\tenex
  \def\bf{\fam\bffam\elevenbf}%
  \def\it{\fam\itfam\elevenit}%
  \textfont\bffam=\elevenbf \scriptfont\bffam=\eightbf
   \scriptscriptfont\bffam=\sixbf
\normalbaselineskip=13.95pt
  \setbox\strutbox=\hbox{\vrule height9.5pt depth4.4pt width0pt\relax}%
  \normalbaselines\rm}

\elevenpoint 

\def\ninepoint{\def\rm{\fam0\ninerm}%
  \textfont0=\ninerm \scriptfont0=\sixrm \scriptscriptfont0=\fiverm
  \textfont1=\ninei \scriptfont1=\sixi \scriptscriptfont1=\fivei
  \textfont2=\ninesy \scriptfont2=\sixsy \scriptscriptfont2=\fivesy
  \textfont3=\tenex \scriptfont3=\tenex \scriptscriptfont3=\tenex
  \def\it{\fam\itfam\nineit}%
  \textfont\itfam=\nineit
  \def\bf{\fam\bffam\ninebf}%
  \textfont\bffam=\ninebf \scriptfont\bffam=\sixbf
   \scriptscriptfont\bffam=\fivebf
\normalbaselineskip=11pt
  \setbox\strutbox=\hbox{\vrule height8pt depth3pt width0pt\relax}%
  \normalbaselines\rm}

\def\eightpoint{\def\rm{\fam0\eightrm}%
  \textfont0=\eightrm \scriptfont0=\sixrm \scriptscriptfont0=\fiverm
  \textfont1=\eighti \scriptfont1=\sixi \scriptscriptfont1=\fivei
  \textfont2=\eightsy \scriptfont2=\sixsy \scriptscriptfont2=\fivesy
  \textfont3=\tenex \scriptfont3=\tenex \scriptscriptfont3=\tenex
  \def\it{\fam\itfam\eightit}%
  \textfont\itfam=\eightit
  \def\bf{\fam\bffam\eightbf}%
  \textfont\bffam=\eightbf \scriptfont\bffam=\sixbf
   \scriptscriptfont\bffam=\fivebf
\normalbaselineskip=12pt
  \setbox\strutbox=\hbox{\vrule height8.5pt depth3.5pt width0pt\relax}%
  \normalbaselines\rm}


\def\eighteenbold{\def\rm{\fam0\eighteenbf}%
  \textfont0=\eighteenbf \scriptfont0=\twelvebf \scriptscriptfont0=\tenbf
  \textfont1=\eighteenmib \scriptfont1=\twelvemib\scriptscriptfont1=\tenmib
  \textfont2=\eighteensyb \scriptfont2=\twelvesyb\scriptscriptfont2=\tensyb
  \textfont3=\eighteenex \scriptfont3=\tenex \scriptscriptfont3=\tenex
  \def\bf{\fam\bffam\eighteenbf}%
  \textfont\bffam=\eighteenbf \scriptfont\bffam=\twelvebf
   \scriptscriptfont\bffam=\tenbf
\normalbaselineskip=20pt
  \setbox\strutbox=\hbox{\vrule height13.5pt depth6.5pt width0pt\relax}%
\everymath {\fam0 }
\everydisplay {\fam0 }
  \normalbaselines\rm}

\def\elevenbold{\def\rm{\fam0\elevenbf}%
  \textfont0=\elevenbf \scriptfont0=\eightbf \scriptscriptfont0=\sixbf
  \textfont1=\elevenmib \scriptfont1=\eightmib \scriptscriptfont1=\sixmib
  \textfont2=\elevensyb \scriptfont2=\eightsyb \scriptscriptfont2=\sixsyb
  \textfont3=\elevenex \scriptfont3=\elevenex \scriptscriptfont3=\elevenex
  \def\bf{\fam\bffam\elevenbf}%
  \textfont\bffam=\elevenbf \scriptfont\bffam=\eightbf
   \scriptscriptfont\bffam=\sixbf
\normalbaselineskip=14pt
  \setbox\strutbox=\hbox{\vrule height10pt depth4pt width0pt\relax}%
\everymath {\fam0 }
\everydisplay {\fam0 }
  \normalbaselines\bf}

\hsize=31pc
\vsize=48pc

\parindent=22pt
\parskip=0pt

\widowpenalty=10000
\clubpenalty=10000

\topskip=12pt

\skip\footins=20pt
\dimen\footins=3in 

\abovedisplayskip=6.95pt plus3.5pt minus 3pt
\belowdisplayskip=\abovedisplayskip


\voffset=7pt\hoffset= .7in

\newif\iftitle

\def\amheadline{\iftitle%
\hbox to\hsize{\hss\currannalsline\hss}\else\line{\ifodd\pageno
\hfill\thetitle\hfill\llap{\elevenrm\folio}\else\rlap{\elevenrm\folio}
\hfill\theauthors\hfill\fi}\fi}

\headline={\amheadline}
\footline={\global\titlefalse}


\def\annalsline#1#2{\vfill\eject
\ifodd\pageno\else 
\line{\hfill}
\vfill\eject\fi
\global\titletrue
\def\currannalsline{\eightrm 
{\eightbf#1} (#2), \thepages}}

\def\titleheadline#1{\def\one{#1}\ifx\one\empty\else
\def\thetitle{{
\let\\ \relax\eightsc\uppercase{#1}}}\fi}

\newif\ifshort

\let\shorttitle\titleheadline

\def\onpages#1#2{\def\thepages{#1--#2}}

\def\thismuchskip[#1]{\vskip#1pt}
\def\ilook{\ifx\next[ \let\go\thismuchskip\else
\let\go\relax\vskip1pt\fi\go}

\def\institution#1{\def\theinstitutions{\vbox{\baselineskip10pt
\def\and{{\eightrm and }}
\def\\{\futurelet\next\ilook}\eightsc #1}}}
\let\institutions\institution

\newwrite\auxfile

\def\startingpage#1{\def\one{#1}\ifx\one\empty\global\pageno=1\else
\global\pageno=#1\fi
\theoremcount=0 \eqcount=0 \sectioncount=0
\openin1 \jobname.aux \ifeof1
\onpages{#1}{???}
\else\closein1 \relax\input \jobname.aux
\onpages{#1}{\lastpage}
\fi\immediate\openout\auxfile=\jobname.aux
}

\def\endarticle{\ifRefsUsed\global\RefsUsedfalse%
\else\vskip21pt\theinstitutions%
\nobreak\vskip8pt
\fi%
\write\auxfile{\string\def\string\lastpage{\the\pageno}}}

\outer\def\bye{\endarticle\par \vfill \supereject \end}

\def\document{\let\fontlist@\relax\let\alloclist@\relax
 \elevenpoint}


\newif\ifacks
\long\def\acknowledgements#1{\def\one{#1}\ifx\one\empty\else
\vskip-\baselineskip%
\global\ackstrue\footnote{\ \unskip}{*#1}\fi}

\def\title#1{\titleheadline{#1}
\vbox to80pt{\vfill
\baselineskip=18pt
\parindent=0pt
\overfullrule=0pt
\hyphenpenalty=10000
\everypar={\hskip\parfillskip\relax}
\hbadness=10000
\def\\ {\vskip1sp}
\eighteenbold#1\vskip1sp}}

\newif\ifauthor

\def\author#1{\vskip11pt
\hbox to\hsize{\hss\tenrm By \tensc#1\ifacks\global\acksfalse*\fi\hss}
\ifshort\else\xdef\theauthors{{\eightsc\uppercase{#1}}}\fi%
\vskip21pt\global\authortrue\everypar={\global\authorfalse\everypar={}}}

\def\twoauthors#1#2{\vskip11pt
\hbox to\hsize{\hss%
\tenrm By \tensc#1 {\tenrm and} #2\ifacks\global\acksfalse*\fi\hss}
\ifshort\else\xdef\theauthors{{\eightsc\uppercase{#1 and #2}}}\fi%
\vskip21pt
\global\authortrue\everypar={\global\authorfalse\everypar={}}}


\newcount\theoremcount
\newcount\sectioncount
\newcount\eqcount

\newif\ifspecialnumon

\def\eqnumber=#1 {\global\eqcount=#1 \global\advance\eqcount by-1\relax}
\def\sectionnumber=#1 {\global\sectioncount=#1
\global\advance\sectioncount by-1\relax}
\def\proclaimnumber=#1 {\global\theoremcount=#1
\global\advance\theoremcount by-1\relax}

\newif\ifsection
\newif\ifsubsection

\def\intro{\global\authorfalse
\centerline{\bf Introduction}\everypar={}\vskip6pt}

\def\elevenboldmath#1{$#1$\egroup}
\def\mathbold{\hbox\bgroup\elevenbold\elevenboldmath}

\def\section#1{\global\theoremcount=0
\global\eqcount=0
\ifauthor\global\authorfalse\else%
\vskip18pt plus 18pt minus 6pt\fi%
{\parindent=0pt
\everypar={\hskip\parfillskip}
\def\\ {\vskip1sp}\elevenpoint\bf%
\ifspecialnumon\global\specialnumonfalse$\rm\spnum$%
\gdef\sectnum{$\rm\spnum$}%
\else\interlinepenalty=10000%
\global\advance\sectioncount by1\relax\the\sectioncount%
\gdef\sectnum{\the\sectioncount}%
\fi. \hskip6pt#1
\vrule width0pt depth12pt}
\hskip\parfillskip
\global\sectiontrue%
\everypar={\global\sectionfalse\global\interlinepenalty=0\everypar={}}%
\ignorespaces

}


\newif\ifspequation

\let\eqno\leqno 

\newif\ifineqalignno
\let\saveleqalignno\leqalignno                        
\def\leqalignno{\let\eqnu\Eeqnu\saveleqalignno}

\let\eqalignno\leqalignno

\def\sectandeqnum{%
\ifspecialnumon\global\specialnumonfalse
$\rm\spnum$\gdef\eqnum{{$\rm\spnum$}}\else\global\firstlettertrue
\global\advance\eqcount by1
\ifappend\applett\else\the\sectioncount\fi.%
\the\eqcount
\xdef\eqnum{\ifappend\applett\else\the\sectioncount\fi.\the\eqcount}\fi}

\def\eqnu{\leqno{\hbox{\elevenrm\ifspequation\else(\fi\sectandeqnum
\ifspequation\global\spequationfalse\else)\fi}}}      

\def\Speqnu{\global\setbox\leqnobox=\hbox{\elevenrm
\ifspequation\else%
(\fi\sectandeqnum\ifspequation\global\spequationfalse\else)\fi}}

\def\Eeqnu{\hbox{\elevenrm
\ifspequation\else%
(\fi\sectandeqnum\ifspequation\global\spequationfalse\else)\fi}}

\newif\iffirstletter
\global\firstlettertrue
\def\eqletter#1{\global\specialnumontrue\iffirstletter\global\firstletterfalse
\global\advance\eqcount by1\fi
\gdef\spnum{\the\sectioncount.\the\eqcount#1}\eqnu}

\newbox\leqnobox
\def\outsideeqnu#1{\global\setbox\leqnobox=\hbox{#1}}

\def\eatone#1{}

\def\dosplit#1#2{\vskip-.5\abovedisplayskip
\setbox0=\hbox{$\let\eqno\outsideeqnu%
\let\eqnu\Speqnu\let\leqno\outsideeqnu#2$}%
\setbox1\vbox{\noindent\hskip\wd\leqnobox\ifdim\wd\leqnobox>0pt\hskip1em\fi%
$\displaystyle#1\mathstrut$\hskip0pt plus1fill\relax
\vskip1pt
\line{\hfill$\let\eqnu\eatone\let\leqno\eatone%
\displaystyle#2\mathstrut$\ifmathqed~~\qed\fi}}%
\copy1
\ifvoid\leqnobox
\else\dimen0=\ht1 \advance\dimen0 by\dp1
\vskip-\dimen0
\vbox to\dimen0{\vfill
\hbox{\unhbox\leqnobox}
\vfill}
\fi}

\everydisplay{\lookforbreak}

\long\def\lookforbreak #1$${\def\mathone{#1}
\expandafter\testforbreak\mathone\splitmath @}

\def\testforbreak#1\splitmath #2@{\def\mathtwo{#2}\ifx\mathtwo\empty%
#1$$%
\ifmathqed\vskip-\belowdisplayskip
\setbox0=\vbox{\let\eqno\relax\let\eqnu\relax$\displaystyle#1$}%
\vskip-\ht0\vskip-3.5pt\hbox to\hsize{\hfill\qed}
\vskip\ht0\vskip3.5pt\fi
\else$$\vskip-\belowdisplayskip
\vbox{\dosplit{#1}{\let\eqno\eatone
\let\splitmath\relax#2}}%
\nobreak\vskip.5\belowdisplayskip
\noindent\ignorespaces\fi}


\newif\ifmathqed



\newcount\linenum
\newcount\colnum

\def\spline{\omit&\multispan{\the\colnum}{\hrulefill}\cr}
\def\colcounter{\ifnum\linenum=1\global\advance\colnum by1\fi}

\def\everyline{\noalign{\global\advance\linenum by1\relax}%
\ifnum\linenum=2\spline\fi}

\def\mtable{\bgroup\offinterlineskip
\everycr={\everyline}\global\linenum=0
\halign\bgroup\vrule height 10pt depth 4pt width0pt
\hfill$##$\hfill\hskip6pt\ifnum\linenum>1
\vrule\fi&&\colcounter\hskip12pt\hfill$##$\hfill\hskip12pt\cr}

\def\endmtable{\crcr\egroup\egroup}




\def\xast{*}
\newcount\intable
\newcount\mathcol
\newcount\savemathcol
\newcount\topmathcol
\newdimen\arrayhspace
\newdimen\arrayvspace

\arrayhspace=8pt 
\arrayvspace=12pt 

\newif\ifdollaron

\def\mathalign#1{\def\arg{#1}\ifx\arg\xast%
\let\go\relax\else\let\go\mathalign%
\global\advance\mathcol by1 %
\global\advance\topmathcol by1 %
\expandafter\def\csname  mathcol\the\mathcol\endcsname{#1}%
\fi\go}

\def\arraypickapart#1]#2*{\if#1c \ifmmode\vcenter\else
\global\dollarontrue$\vcenter\fi\else%
\if#1t\vtop\else\if#1b\vbox\fi\fi\fi\bgroup%
\def\one{#2}}

\def\arraystrut{\vrule height .7\arrayvspace depth .3\arrayvspace width 0pt}

\def\array#1#2*{\def\firstarg{#1}%
\if\firstarg[ \def\two{#2} \expandafter\arraypickapart\two*\else%
\ifmmode\vcenter\else\vbox\fi\bgroup \def\one{#1#2}\fi%
\global\everycr={\noalign{\global\mathcol=\savemathcol\relax}}%
\def\\ {\cr}%
\global\advance\intable by1 %
\ifnum\intable=1 \global\mathcol=0 \savemathcol=0 %
\else \global\advance\mathcol by1 \savemathcol=\mathcol\fi%
\expandafter\mathalign\one*%
\mathcol=\savemathcol %
\halign\bgroup&\hskip.5\arrayhspace\arraystrut%
\global\advance\mathcol by1 \relax%
\expandafter\if\csname mathcol\the\mathcol\endcsname r\hfill\else%
\expandafter\if\csname mathcol\the\mathcol\endcsname c\hfill\fi\fi%
$\displaystyle##$%
\expandafter\if\csname mathcol\the\mathcol\endcsname r\else\hfill\fi\relax%
\hskip.5\arrayhspace\cr}

\def\endarray{\crcr\egroup\egroup%
\global\mathcol=\savemathcol %
\global\advance\intable by -1\relax%
\ifnum\intable=0 %
\ifdollaron\global\dollaronfalse $\fi
\loop\ifnum\topmathcol>0 %
\expandafter\def\csname  mathcol\the\topmathcol\endcsname{}%
\global\advance\topmathcol by-1 \repeat%
\global\everycr={}\fi%
}

\def\big#1{{\hbox{$\left#1\vbox to 10pt{}\right.\n@space$}}}
\def\Big#1{{\hbox{$\left#1\vbox to 13pt{}\right.\n@space$}}}
\def\bigg#1{{\hbox{$\left#1\vbox to 16pt{}\right.\n@space$}}}
\def\Bigg#1{{\hbox{$\left#1\vbox to 19pt{}\right.\n@space$}}}


\def\figcaption#1#2#3{\topinsert
\vskip4pt 
\vbox to#3{\vfill}\vskip1sp
\setbox0=\hbox{\eightsc Figure #1.\hskip12pt\eightpoint #2}
\ifdim\wd0>\hsize
\noindent\eightsc Figure #1.\hskip12pt\eightpoint #2
\else
\centerline{\eightsc Figure #1.\hskip12pt\eightpoint #2}
\fi
\vskip16pt
\endinsert}

\def\wfig#1#2#3{\topinsert
\vskip4pt 
\hbox to\hsize{\hss\vbox{\hrule height .25pt width #3
\hbox to #3{\vrule width .25pt height #2\hfill\vrule width .25pt height #2}
\hrule height.25pt}\hss}
\vskip1sp
\centerline{\eightsc Figure #1}
\vskip16pt
\endinsert}

\def\wfigcaption#1#2#3#4{\topinsert
\vskip4pt 
\hbox to\hsize{\hss\vbox{\hrule height .25pt width #4
\hbox to #4{\vrule width .25pt height #3\hfill\vrule width .25pt height #3}
\hrule height.25pt}\hss}
\vskip1sp
\setbox0=\hbox{\eightsc Figure #1.\hskip12pt\eightpoint\rm #2}
\ifdim\wd0>\hsize
\noindent\eightsc Figure #1.\hskip12pt\eightpoint\rm #2\else
\centerline{\eightsc Figure #1.\hskip12pt\eightpoint\rm #2}\fi
\vskip16pt
\endinsert}

\def\tabcaption#1#2{\vskip6pt
\setbox0=\hbox{\eightsc Table #1.\hskip12pt\eightpoint #2}
\ifdim\wd0>\hsize
\noindent\eightsc Table #1.\hskip12pt\eightpoint #2
\else
\centerline{\eightsc Table #1.\hskip12pt\eightpoint #2}
\fi
\vskip6pt}

\def\endinsert{\egroup\if@mid\dimen@\ht\z@\advance\dimen@\dp\z@
\advance\dimen@ 12\p@\advance\dimen@\pagetotal\ifdim\dimen@ >\pagegoal
\@midfalse\p@gefalse\fi\fi\if@mid\smallskip\box\z@\bigbreak\else
\insert\topins{\penalty 100 \splittopskip\z@skip\splitmaxdepth\maxdimen
\floatingpenalty\z@\ifp@ge\dimen@\dp\z@\vbox to\vsize {\unvbox \z@
\kern -\dimen@ }\else\box\z@\nobreak\smallskip\fi}\fi\endgroup}

\def\pagecontents{
\ifvoid\topins \else\iftitle\else
\unvbox \topins \fi\fi \dimen@ =\dp \@cclv \unvbox
\@cclv
\ifvoid\topins\else\iftitle\unvbox\topins\fi\fi
\ifvoid \footins \else \vskip \skip \footins \footnoterule
\unvbox \footins \fi \ifr@ggedbottom \kern -\dimen@ \vfil \fi}


\newif\ifappend

\def\appendix#1#2{\def\applett{#1}\def\two{#2}%
\global\appendtrue
\global\theoremcount=0
\global\eqcount=0
\vskip18pt plus 18pt
\vbox{\parindent=0pt
\everypar={\hskip\parfillskip}
\def\\ {\vskip1sp}\elevenbold Appendix%
\ifx\applett\empty\gdef\applett{A}\ifx\two\empty\else.\fi%
\else\ #1.\fi\hskip6pt#2\vskip12pt}%
\global\sectiontrue%
\everypar={\global\sectionfalse\everypar={}}\nobreak\ignorespaces}

\newif\ifRefsUsed
\long\def\references{\global\RefsUsedtrue\vskip21pt
\theinstitutions
\global\everypar={}\global\bibnum=0
\vskip20pt\goodbreak\bgroup
\vbox{\centerline{\eightsc References}\vskip6pt}%
\ifdim\maxbibwidth>0pt
\leftskip=\maxbibwidth%
\parindent=-\maxbibwidth%
\else
\leftskip=18pt%
\parindent=-18pt%
\fi
\ninepoint
\frenchspacing
\nobreak\ignorespaces\everypar={\amref}%
}

\def\endreferences{\vskip1sp\egroup\global\everypar={}%
\nobreak\vskip8pt\vbox{\thereceived\therevised}
}

\newcount\bibnum

\def\amref#1 {\global\advance\bibnum by1%
\immediate\write\auxfile{\string\expandafter\string\def\string\csname
\space #1croref\string\endcsname{[\the\bibnum]}}%
\leavevmode\hbox to18pt{\hbox to13.2pt{\hss[\the\bibnum]}\hfill}}

\def\bibline{\hbox to30pt{\hrulefill}\/\/}

\def\name#1{{\eightsc#1}}

\newdimen\maxbibwidth
\def\AuthorRefNames [#1] {%
\immediate\write\auxfile{\string\def\string\cite\string##1{[\string##1]}}

\def\amref{\spamref}
\setbox0=\hbox{[#1] }\global\maxbibwidth=\wd0\relax}

\def\spamref[#1] {\leavevmode\hbox to\maxbibwidth{\hss[#1]\hfill}}


\def\footnoterule{\kern-3pt\hrule width1in height.5pt\kern2.5pt}

\def\footnote#1#2{%
\plainfootnote{#1}{{\eightpoint\normalbaselineskip11pt
\normalbaselines#2}}}

\def\vfootnote#1{%
\insert \footins \bgroup \eightpoint\baselineskip11pt
\interlinepenalty \interfootnotelinepenalty
\splittopskip \ht \strutbox \splitmaxdepth \dp \strutbox \floatingpenalty
\@MM \leftskip \z@skip \rightskip \z@skip \spaceskip \z@skip
\xspaceskip \z@skip
{#1}$\,$\footstrut \futurelet \next \fo@t}


\newif\iffirstadded
\newif\ifadded

\def\addedlett{}

\def\alltheoremnums{%
\ifspecialnumon\global\specialnumonfalse
\ifadded\global\addedfalse
\iffirstadded\global\firstaddedfalse
\global\advance\theoremcount by1 \fi
\ifappend\applett\else\the\sectioncount\fi.\the\theoremcount\addedlett%
\xdef\theoremnum{\ifappend\applett\else\the\sectioncount\fi.%
\the\theoremcount\addedlett}%
\else$\rm\spnum$\def\theoremnum{{$\rm\spnum$}}\fi%
\else\global\firstaddedtrue
\global\advance\theoremcount by1
\ifappend\applett\else\the\sectioncount\fi.\the\theoremcount%
\xdef\theoremnum{\ifappend\applett\else\the\sectioncount\fi.%
\the\theoremcount}\fi}

\def\allcorolnums{%
\ifspecialnumon\global\specialnumonfalse
\ifadded\global\addedfalse
\iffirstadded\global\firstaddedfalse
\global\advance\corolcount by1 \fi
\the\corolcount\addedlett%
\else$\rm\spnum$\def\corolnum{$\rm\spnum$}\fi%
\else\global\advance\corolcount by1
\the\corolcount\fi}


\newcount\corolcount
\def\xcorol{Corollary}
\def\xtheorem{Theorem}
\def\xmaintheorem{Main Theorem}

\newif\ifthtitle

\let\saverparen)
\let\savelparen(
\def\rmparenl{{\rm(}}
\def\rmparenr{{\rm\/)}}
{
\catcode`(=13
\catcode`)=13
\gdef\makeparensRM{\catcode`(=13\catcode`)=13\let(=\rmparenl%
\let)=\rmparenr%
\everymath{\let(\savelparen%
\let)\saverparen}%
\everydisplay{\let(\savelparen%
\let)\saverparen\lookforbreak}}}

\medskipamount=8pt plus.1\baselineskip minus.05\baselineskip

\def\rmtext#1{\hbox{\rm#1}}

\def\proclaim#1{\vskip-\lastskip
\def\one{#1}\ifx\one\xtheorem\global\corolcount=0\fi
\ifsection\global\sectionfalse\vskip-6pt\fi
\medskip
{\elevensc#1}%
\ifx\one\xmaintheorem\global\corolcount=0
\gdef\theoremnum{Main Theorem}\else%
\ifx\one\xcorol\ \allcorolnums\else\ \alltheoremnums\fi\fi%
\ifthtitle\ \global\thtitlefalse{\rm(\thethtitle)}\fi.%
\hskip1em\bgroup\let\text\rmtext\makeparensRM\it\ignorespaces}

\def\nonumproclaim#1{\vskip-\lastskip
\def\one{#1}\ifx\one\xtheorem\global\corolcount=0\fi
\ifsection\global\sectionfalse\vskip-6pt\fi
\medskip
{\elevensc#1}.\ifx\one\xmaintheorem\global\corolcount=0
\gdef\theoremnum{Main Theorem}\fi\hskip.5pc%
\bgroup\it\makeparensRM\ignorespaces}

\def\endproclaim{\egroup\medskip}


\def\xproof{Proof}
\def\xremark{Remark}
\def\xcase{Case}
\def\xsubcase{Subcase}
\def\xconjecture{Conjecture}
\def\xstep{Step}
\def\xof{of}

\def\deconstruct#1 #2 #3 #4 #5 @{\def\one{#1}\def\two{#2}\def\three{#3}%
\def\four{#4}%
\ifx\two\empty #1\else%
\ifx\one\xproof%
\ifx\two\xof%
  \ifx\three\xcorol Proof of Corollary \rm#4\else%
     \ifx\three\xtheorem Proof of Theorem \rm#4\else\xone\fi%
  \fi\fi%
\else\xone\fi\fi.}

\def\pickup#1 {\def\this{#1}%
\ifx\this\xproof\global\let\go\demoproof
\global\let\enddemo\endproof\else
\ifx\this\xremark\global\let\go\demoremark\else
\ifx\this\xcase\global\let\go\demostep\else
\ifx\this\xsubcase\global\let\go\demostep\else
\ifx\this\xconjecture\global\let\go\demostep\else
\ifx\this\xstep\global\let\go\demostep\else
\global\let\go\demoproof\fi\fi\fi\fi\fi\fi}

\newif\ifnonum
\def\demo#1{\vskip-\lastskip
\ifsection\global\sectionfalse\vskip-6pt\fi
\def\one{#1 }\def\two{#1*}%
\setbox0=\hbox{\expandafter\pickup\one}\expandafter\go\two}

\def\numbereddemo#1{\vskip-\lastskip
\ifsection\global\sectionfalse\vskip-6pt\fi
\def\two{#1*}%
\expandafter\demoremark\two}

\def\demoproof#1*{\medskip\def\xone{#1}
{\ignorespaces\it\expandafter\deconstruct\xone {} {} {} {} {} @%
\unskip\hskip6pt}\rm\ignorespaces}

\def\demoremark#1*{\medskip
{\it\ignorespaces#1\/} \ifnonum\global\nonumtrue\else
 \alltheoremnums\unskip.\fi\hskip1pc\rm\ignorespaces}

\def\demostep#1 #2*{\vskip4pt
{\it\ignorespaces#1\/} #2.\hskip1pc\rm\ignorespaces}

\def\enddemo{\medskip}

\def\endproof{\ifmathqed\global\mathqedfalse\medskip\else
\parfillskip=0pt~~\hfill\qed\medskip
\fi\global\parfillskip0pt plus 1fil\relax
\gdef\enddemo{\medskip}}

\def\qed{\vbox{\hrule\hbox{\vrule height6pt\hskip6pt\vrule}\hrule}}


\def\proofbox{\parfillskip=0pt~~\hfill\qed\vskip1sp\parfillskip=
0pt plus 1fil\relax}







\def\stripbs#1#2*{\def\one{#2}}

\def\emptyspace{ }
\def\nextthing{}
\def\newline{***}
\def\eatone#1{ }

\def\lookatspace#1{\ifcat\noexpand#1\ \else%
\gdef\nextthing{}\xdef\next{#1}%
\ifx\next\emptyspace%
\let\nextthing\emptyspace\else\ifx\next\newline%
\gdef\nextthing{\eatone}\fi\fi\fi\egroup\nextthing#1}

{\catcode`\^^M=\active%
\gdef\spacer{\bgroup\catcode`\^^M=\active%
\let^^M=\newline\obeyspaces\lookatspace}}

\def\ref#1{\seeifdefined{#1}\expandafter\csname\one\endcsname\spacer}

\def\cite#1{\expandafter\ifx\csname#1croref\endcsname\relax[??]\else
\csname#1croref\endcsname\fi\spacer}


\def\seeifdefined#1{\expandafter\stripbs\string#1croref*%
\crorefdefining{#1}}

\newif\ifcromessage
\global\cromessagetrue

\def\crorefdefining#1{\ifdefined{\one}{}
{\ifcromessage\global\cromessagefalse%
\message{\spaces\spaces\spaces\spaces\spaces\spaces\spaces}%
\message{<Undefined reference.}%
\message{Please TeX file once more to have accurate cross-references.>}%
\message{\spaces\spaces\spaces\spaces\spaces\spaces\spaces}\fi[??]}}

\def\label#1#2*{\gdef\ctest{#2}%
\xdef\currlabel{\string#1croref}
\expandafter\seeifdefined{#1}%
\ifx\empty\ctest%
\xdef\labelnow{\write\auxfile{\noexpand\def\currlabel{\the\pageno}}}%
\else\xdef\labelnow{\write\auxfile{\noexpand\def\currlabel{#2}}}\fi%
\labelnow}

\def\ifdefined#1#2#3{\expandafter\ifx\csname#1\endcsname\relax%
#3\else#2\fi}




\def\articlecontents{
\vskip20pt\centerline{\bf Table of Contents}\everypar={}\vskip6pt
\bgroup \leftskip=3pc \parindent=-2pc
\def\item##1{\vskip1sp\indent\hbox to2pc{##1.\hfill}}}

\def\endcontents{\vskip1sp\leftskip=0pt\egroup}

\def\journalcontents{\vfill\eject
\def\currannalsline{\hfill}
\global\titletrue
\vglue3.5pc
\centerline{\tensc\hskip12pt TABLE OF CONTENTS}\everypar={}\vskip30pt
\bgroup \leftskip=34pt \rightskip=-12pt \parindent=-22pt
  \def\\ {\vskip1sp\noindent}
\def\pagenum##1{\unskip\parfillskip=0pt\dotfill##1\vskip1sp
\parfillskip=0pt plus 1fil\relax}
\def\name##1{{\tensc##1}}}


\institution{}
\onpages{0}{0}
\def\lastpage{???}
\def\thetitle{Title ???}
\def\theauthors{Authors ???}
\def\thereceived{}
\def\therevised{}

\gdef\split{\relaxnext@\ifinany@\let\next\insplit@\else
 \ifmmode\ifinner\def\next{\onlydmatherr@\split}\else
 \let\next\outsplit@\fi\else
 \def\next{\onlydmatherr@\split}\fi\fi\let\eqnu\xspliteqnu\next}

\gdef\align{\relaxnext@\ifingather@\let\next\galign@\else
 \ifmmode\ifinner\def\next{\onlydmatherr@\align}\else
 \let\next\align@\fi\else
 \def\next{\onlydmatherr@\align}\fi\fi\let\eqnu\xspliteqnu\next}

\def\spliteqnu{{\tenrm\sectandeqnum}\relax}

\def\xspliteqnu{\tag\spliteqnu}

\catcode`@=12

\document

\annalsline{May}{1995}
\startingpage{1}     

\comment
\nopagenumbers
\headline{\ifnum\pageno=1\hfil\else \rightheadline\fi}
\def\rightheadline{\hfil\eightit
The Macdonald conjecture
\quad\eightrm\folio}

\voffset=2\baselineskip
\endcomment


%
%
%
%
%

\def\for{\  \hbox{ for } \ }
\def\if{ \ \hbox{ if } \ }
\def\when{ \ \hbox{ when } \ }
\def\where{\  \hbox{ where } \ }
\def\and{\  \hbox{ and } \ }

\def\equal{\buildrel  def \over =}

\def\la{\lambda}
\def\La{\Lambda}
\def\om{\omega}
\def\Om{\Omega}

\def\th{\theta}
\def\al{\alpha}
\def\be{\beta}
\def\ga{\gamma}
\def\ep{\epsilon}

\def\de{\delta}

\def\si{\sigma}
\def\Si{\Sigma}
\def\Ga{\Gamma}
\def\ze{\zeta}

\def\pa{\partial}

\def\vph{\varphi}

\def\vep{\varepsilon}

\def\tal{\tilde{\alpha}}
\def\tbe{\tilde{\beta}}
\def\tde{\tilde{\delta}}
\def\tpi{\tilde{\pi}}
\def\tPi{\tilde{\Pi}}
\def\tV{\tilde{V}}

\def\tw{\tilde w}

\def\tB{\tilde B}

\def\tV{\tilde V}
\def\tz{\tilde z}
\def\tb{\tilde b}

\def\hH{\hat{H}}

\def\hY{\hat{Y}}

\def\hT{\hat{T}}

\def\hw{\hat{w}}

\def\hv{\hat{v}}

\def\C{\bold{C}}
\def\Q{\bold{Q}}
\def\B{\bold{B}}
\def\R{\bold{R}}
\def\N{\bold{N}}
\def\Z{\bold{Z}}

\def\one{\bold{1}}

\def\0{\bold{0}}

\def\C{\hbox{\bf C}}


\def\f{\Cal{F}}
\def\t{\Cal{T}}

\def\l{\Cal{L}}
\def\m{\Cal{M}}

\def\d{\Cal{D}}
\def\p{\Cal{P}}
\def\a{\Cal{A}}

\def\y{\Cal{Y}}
\def\e{\Cal{E}}
\def\v{\Cal{V}}

\def\x{\Cal{X}}
\def\g{\Cal{G}}

\def\w{\Cal{W}}

\font\germ=eufb10 
\def\goth#1{\hbox{\germ #1}}

\def\TT{\goth{T}}
\def\HH{\goth{H}}

\def\BB{\goth{B}}
\def\AA{\goth{A}}

\font\smm=msbm10 at 12pt
\def\symbol#1{\hbox{\smm #1}}
\def\lsmash{{\symbol n}}



\title
{Macdonald's Evaluation Conjectures,\\
 Difference Fourier Transform,\\
and applications}

\shorttitle{ Evaluation Conjectures}

\acknowledgements{
Partially supported by NSF grant DMS--9301114}

\author{ Ivan Cherednik}

\institutions{
Math. Dept, University of North Carolina at Chapel Hill,
 N.C. 27599-3250
\\ Internet: chered\@math.unc.edu
}


\intro 
%
%
%
%
%
\vfil
 Generalizing
the characters of compact simple Lie groups,
Ian Macdonald introduced in [M1,M2] and other works
remarkable orthogonal symmetric
polynomials dependent on the parameters
$q,t$. He came up with three main conjectures  formulated
for arbitrary root systems. A new approach to
the Macdonald theory was suggested in [C1] on the basis of
(double) affine Hecke algebras and related difference operators.
In [C2] the norm conjecture
(including the celebrated constant term conjecture [M3])
was proved for all (reduced) root systems.
This paper contains the proof of
the remaining two (the duality and evaluation conjectures),
the recurrence relations, and basic results on Macdonald's
polynomials at roots of unity.
In the next paper the same questions will be considered for
the non-symmetric polynomials.

The evaluation conjecture (now a theorem) is in fact
a $q,t$-generalization
of the classic Weyl dimension formula. One can expect
interesting applications of this theorem since
the so-called $q$-dimensions are undoubtedly
important.  To demonstrate  deep relations to
the representation theory we prove the
Recurrence Theorem  connected with
decomposing of the tensor products of represenations
of compact Lie groups in terms of irreducible ones.
The arising $q,t$-multiplicities
 are in fact the
coefficients of our difference
operators (one needs the duality to establish
this).
It is  likely that we can incorporate
the Kac-Moody case as well.  The necessary technique
was developed in [C4].

The  duality theorem (in its complete form) states
that the  difference
zonal $q,t$-Fourier transform
is self-dual (its reproducing kernel is symmetric).
In this paper we introduce the transform formally
in terms of double affine
Hecke algebras. The  self-duality is directly related to
the interpretation of these algebras via the so-called
elliptic braid groups (the Fourier involution  turns into the
transposition of the periods of an elliptic curve).
It is not very surprising since this interpretation is actually
the monodromy representation of
the double affine (elliptic) Knizhnik-Zamolodchikov equation
from [C6].

\vskip 5pt

The classical zonal (spherical) Fourier transform
plays one of the main roles in the harmonic analysis
on symmetric spaces $G/K$. It sends the radial  parts of
$K$-invariant differential operators to the corresponding symmetric
 polynomials (the Harish-Chandra isomorphism)
 and is not self-dual. The calculation
of its inverse (the Plancherel theorem) is
 important and involved.
This transform is the limiting case of our construction
when $q\to 1$ and $t=q^k$  for certain special
$k$.

Considering the rational theory (replacing the trigonometric
coefficients of the operators by their rational degenerations),
Charles Dunkl introduced
the generalized Hankel transform which appeared to be
 self-dual [D,J]. Geometrically, it is the case of
the tangent space $T_e(G/K)$ with the adjoint action of $K$,
taken instead of $G/K$ (see [H]).
We demonstrate in this paper that one can save this very
important property in the main (trigonometric) setup
going to
 the difference counterparts of the zonal Laplace operators.
At the moment, it is mostly an algebraic observation. The
difference-analytical aspects were  not studied systematically
(see Section 4  where we discuss Schwartz functions).

\vskip 5pt

As to the differential theory, the evaluation conjecture
 was proved  by Eric Opdam
[O1] (see also [O2]) together with other
 Macdo\-nald- Mehta
conjectures excluding
the duality con\-jec\-ture which
collapses.
 He used
the Heckman-Opdam operators
 (including the shift operator - see [O1,He]).
We use their generalizations from [C1,C2]
defined by means of double
affine Hecke algebras. We mention that the proof of the norm
conjecture from  [C2] was based mainly on
the  relations and properties of affine
Hecke algebras. Only
in this paper  double Hecke algebras  work at their full potential
 ensuring the duality.
\vfil

Concerning the  open questions in the Macdonald theory,
 we will say a little bit about
 $A_n$. First of all, the Macdonald
polynomials (properly normalized) are $q\leftrightarrow t$
symmetric and satisfy various additional
relations. Then there are
very interesting positivity conjectures
(Macdonald [M1], Garsia, Haiman [GH]).
Moreover they can be interpreted as generalized characters
(Etingof, Kirillov [EK1]).
In the differential setting, these
 polynomials (Jack polynomials)
are also quite remarkable
(Stanley, Hanlon).
 By the way, due to Andrews
one can add $n$ new $t$-parameters
and still the
constant term conjecture will hold (Bressoud and
Zeilberger),  but at the moment we have no definition of the
associated
 orthogonal polynomials.
Hopefully some of these properties can be extended to
arbitrary roots.

\vfil
Recently Alexander Kirillov, Jr. established that in the
case of $A_{n}$ the action
of $SL_2(\Z)$ appearing in the theory of quantum $SL_{n+1}$ at
roots of unity (see [MS]) leads to the projective representations of this
groop expressed in terms of special values of the Macdonald
polynomials (when $q$ is a root of unity, $k\in \Z_+$).
We demonstrate   that this result can
be naturally rediscovered in the frameworks of the
present paper and generalized to arbitrary root systems.
 When $q=t$ the Macdonald polynomials are the
characters (Schur functions for $A_n$)
and we arrive at Theorem 13.8 from [K].
 Kirillov's work [Ki] is expected to be directly
connected with the equivalence of the quantum groups and
Kac-Moody algebras due to Kazhdan - Lusztig [KL2] (in the
case of $A_n$). In our setting
it should result from the calculation of
the monodromy representation of
the double affine  Knizhnik-Zamolodchikov equation
(for arbitrary root systems).

We note that this paper is a part of a new program
in the harmonic analysis of
symmetric spaces based  on certain
remarkable representations of  Hecke algebras
in terms of Dunkl and
Demazure operators instead of Lie groups and Lie algebras.
It gave already a $k$-parametric deformation of
the classical  theory (see
[O1,He,C5]) directly connected with the so-called
quantum many-body problem
(Calogero,  Sutherland, Moser, Olshanetsky,
 Perelomov). Then it was extended
(in the algebraic context)
 to the difference, elliptic,
and finally to
the difference-elliptic case [C4]
presumably corresponding to the quantum Kac-Moody algebras.
We say presumably because the harmonic analysis
for the latter algebras does not exist.


{\it The duality-evaluation conjecture.}
Let $R=\{\al\}   \subset \R^n$ be a root system of type $A,B,...,F,G$
with respect to a euclidean form $(z,z')$ on $\R^n \ni z,z'$,
$W$ the Weyl group  generated by the the reflections $s_\al$.
We assume that $(\al,\al)=2$ for long $\al$.
Let us  fix the set $R_{+}$ of positive  roots ($R_-=-R_+$),
the corresponding simple
roots $\al_1,...,\al_n$, and  their dual counterparts
$a_1 ,..., a_n,  a_i =\al_i^\vee, \where \al^\vee =2\al/(\al,\al)$.
The dual fundamental weights
$b_1,...,b_n$  are determined from the relations  $ (b_i,\al_j)=
\de_i^j $ for the
Kronecker delta. We will also introduce the dual root system
$R^\vee =\{\al^\vee, \al\in R\}, R^\vee_+$, and the lattices
$$
\eqalignno{
& A=\oplus^n_{i=1}\Z a_i \subset B=\oplus^n_{i=1}\Z b_i,
}
$$
  $A_\pm, B_\pm$  for $\Z_{\pm}=\{m\in\Z, \pm m\ge 0\}$
instead of $\Z$. (In the standard notations, $A= Q^\vee,\
B = P^\vee $ - see [B].)  Later on,
$$
\eqalign{
&\nu_{\al}=\nu_{\al^\vee}\ =\ (\al,\al),\  \nu_i\ =\ \nu_{\al_i}, \
\nu_R\ = \{\nu_{\al}, \al\in R\}, \cr
&\rho_\nu\ =\ (1/2)\sum_{\nu_{\al}=\nu} \al \ =
\ (\nu/2)\sum_{\nu_i=\nu}  b_i, \for\al\in R_+,\cr
&r_\nu\ =\ \rho_\nu^\vee \ =\ (2/\nu)\rho_\nu\ =\
\sum_{\nu_i=\nu}  b_i,\quad 2/\nu=1,2,3.
}
\eqnu
$$

Let us put formally $x_i=exp({b_i}),\  x_b=exp(b)= \prod_{i=1}^n
x_i^{k_i} \for b=\sum_{i=1}^n k_i b_i$,  and introduce the algebra
$\C(q,t)[x]$  of polynomials in terms of $x_i^{\pm 1}$ with the
coefficients belonging to the field $\C(q,t)$ of rational functions
in terms  of indefinite complex parameters $q, t_\nu,
\nu\in \nu_R$ (we will put $t_\al=t_{\nu_\al}=t_{\al^\vee}$).
The coefficient of $x^0=1$ ({\it the constant term})
will be denoted by $\langle \  \rangle$. The following product is a
Laurent series in $x$ with the coefficients in  the algebra
$\C[t][[q]]$ of formal series in $q$ over polynomials in $t$:
$$
\eqalign{
&\mu\ =\ \prod_{a \in R_+^\vee}
\prod_{i=0}^\infty {(1-x_aq_a^{i}) (1-x_a^{-1}q_a^{i+1})
\over
(1-x_a t_aq_a^{i}) (1-x_a^{-1}t_a^{}q_a^{i+1})},
}
\eqnu
\label\mu\eqnum*
$$
where $q_a=q_{\nu}=q^{2/\nu} \for \nu=\nu_a$.
We note that  $\mu\in
\C(q,t)[x]$ if $t_\nu=q_\nu^{k_\nu}$ for $k_\nu\in \Z_+$.
The coefficients of $\mu_1\equal \mu/\langle \mu \rangle$
are from $\C(q,t)$, where the formula for the
constant term of $\mu$ is as follows
(see [C2]):
$$
\eqalign{
&\langle\mu\rangle\ =\ \prod_{a \in R_+^\vee}
\prod_{i=1}^\infty {(1-x_a(t^\rho)q_a^{i})^2
\over
(1-x_a(t^\rho) t_aq_a^{i}) (1-x_a(t^\rho) t_a^{-1}q_a^{i})}.
}
\eqnu
\label\consterm\eqnum*
$$
Here and further
 $x_{b}(t^{\pm\rho}q^{c})=
q^{(b,c)}\prod_\nu t_\nu^{\pm(b,\rho_\nu)}$. We note that
$\mu_1^*\ =\ \mu_1$  with respect to the involution
$$
 x_b^*\ =\  x_{-b},\ t^*\ =\ t^{-1},\ q^*\ =\ q^{-1}.
$$

The
{\it monomial symmetric functions}
$m_{b}\ =\ \sum_{c\in W(b)}x_{c}$ for $b\in B_-$
form a base of the space
 $\C[x]^W$ of all $W$-invariant polynomials.
Setting ,
$$
\eqalignno{
&\langle f,g\rangle\ =\langle \mu_1 f\ {g}^*\rangle\ =\
\langle g,f\rangle^* \for
f,g \in \C(q,t)[x]^W,
&\eqnu
\label\innerpro\eqnum*
}
$$
 we  introduce the {\it Macdonald
polynomials} $p_b(x),\   b \in B_-$, by means of
the conditions
$$
\eqalignno{
&p_b-m_b\ \in\ \oplus_c\C(q,t)m_{c},\
\langle p_b, m_{c}\rangle = 0, \for c\succ b  &\eqnu
 \cr
&\hbox{where\ }
 c\in B_-,\
 c\succ b \hbox{\ means\ that\ }  c-b \in A_+, c\neq b.
\label\macd\eqnum*
}
$$
They can be determined by the Gram - Schmidt process
because the (skew Macdonald) pairing (see [M1,M2,C2])
is non-degenerate
 and form a
basis in $\C(q,t)[x]^W$. As it was established by Macdonald
they are pairwise orthogonal for arbitrary $b\in B_-$.
We note that $p_b$ are ''real'' with respect to the formal
 conjugation sending $q\to q^{-1},\ t\to t^{-1}$.
It makes our definition compatible with  Macdonald's
original one (his $\mu$ is somewhat different).

\proclaim{Main Theorem}
Given $b,c\in B_-$ and the corresponding Macdonald
polynomials $p_b, p_c$,
$$
\eqalignno{
&p_b(t^{-\rho} q^{c})p_c(t^{-\rho})\ =\
p_c(t^{-\rho} q^{b})p_b(t^{-\rho}),
&\eqnu\cr
\label\PP\eqnum*
&p_b(t^{-\rho}) = p_b(t^{\rho}) =
x_b(t^{\rho})\prod_{a\in R_+^\vee}
\prod_{ j=1}^{-(a^\vee,b)}
\Bigl(
{
1- q_a^{j-1}t_a x_a(t^\rho)
 \over
1- q_a^{j-1}x_a(t^\rho)
}
\Bigr) =\cr
&x_b(t^\rho)
\prod_{a\in R_+^\vee, 0\le j< \infty}
\Bigl(
{
(1-q_a^{j}x_a(t^\rho q^{-b}))
(1-t_aq_a^{j}x_a(t^\rho))
 \over
(1-t_aq_a^{j}x_a(t^\rho q^{-b}))
(1- q_a^{j}x_a(t^\rho))
}
\Bigr).&\eqnu\cr
\label\EV\eqnum*
}
$$
\label\MAIN\theoremnum*
\endproclaim
The right hand side of (\ref\EV) is a rational function
in terms of $ q,t$ ($a^\vee=2a/(a,a)\in R$).
 We mention that there is a straightforward passage
 to non-reduced root systems,  the non-symmetric Macdonald
polynomials, and
to  $\mu$  introduced for
$\al\in R_+$ instead of $a\in R_+^\vee$ (see [C2]). As to the latter
case, it is necessary just to replace the indices $a$ by $\al$
($ q_a \to q, \rho\to r$) in all final formulas for $\{p_a\}$.

The second formula was conjectured by Macdonald (see (12.10),[M2]).
 He also formulated an equivalent version of (\ref\PP)
in one of his lectures (1991). Both statements were
 established for $A_n$  by Koornwinder in 1988 (his proof
was not published) and by Macdonald (to be published).
Recently the paper by Etingof and Kirillov [EK2]
appeared  were they use their interpretation of the Macdonald
polynomials to check the above theorem (and the norm conjecture)
in the case of $A_n$. As to other root systems,
it seems that almost
nothing was
known (excluding $BC_1$ and certain
 special values of the
parameters).

Concerning the notations, we use $t$ for Hecke algebras
and $q$ for difference operators. So in this paper
(in contrast to [C2]) we switch to the
standard meaning of these letters in the papers on
$q$-orthogonal polynomials and $q$-functions.

The author thanks G. Heckman for useful remarks
and A. Kirillov, Jr. for  discussion and the opportunity
to know about his results before their publication.
The paper was completed at Harvard University.
I am grateful for the kind invitation and hospitality.
I acknowledge
my special indebtedness to D. Kazhdan who inspired this
paper a lot.

%
%
%
%
\section { Double affine Hecke algebras}
The vectors $\ \tal=[\al,k] \in
\R^n\times \R \subset \R^{n+1}$
for $\al \in R, k \in \Z $
form the {\it affine root system}
$R^a \supset R$ ( $z\in \R^n$ are identified with $ [z,0]$).
We add  $\al_0 \equal [-\th,1]$ to the  simple roots
for the {\it maximal root} $\th \in R$.
The corresponding set $R^a_+$ of positive roots coincides
with $R_+\cup \{[\al,k],\  \al\in R, \  k > 0\}$.

We denote the Dynkin diagram and its affine completion with
$\{\al_j,0 \le j \le n\}$ as the vertices by $\Ga$ and $\Ga^a$.
Let $m_{ij}=2,3,4,6$\  if $\al_i\and\al_j$ are joined by 0,1,2,3 laces
respectively.
The set of
the indices of the images of $\al_0$ by all
the automorphisms of $\Ga^a$ will be denoted by $O$ ($O=\{0\}
\for E_8,F_4,G_2$). Let $O^*={r\in O, r\neq 0}$.
The elements $b_r$ for $r\in O^*$ are the so-called minuscule
weights ($(b_r,\al)\le 1$ for
$\al \in R_+$).

Given $\tal=[\al,k]\in R^a,  \ b \in B$, let
$$
\eqalignno{
&s_{\tal}(\tz)\ =\  \tz-(z,\al^\vee)\tal,\
\ b'(\tz)\ =\ [z,\ze-(z,b)]
&\eqnu
}
$$
for $\tz=[z,\ze] \in \R^{n+1}$.

The {\it affine Weyl group} $W^a$ is generated by all $s_{\tal}$
(we write $W^a = <s_{\tal}, \tal\in R_+^a>)$. One can take
the simple reflections $s_j=s_{\al_j}, 0 \le j \le n,$ as its
generators and introduce the corresponding notion of the
length. This group is
the semi-direct product $W\lsmash A'$ of
its subgroups $W=<s_\al,
\al \in R_+>$ and $A'=\{a', a\in A\}$, where
$$
\eqalignno{
& a'=\ s_{\al}s_{[\al,1]}=\ s_{[-\al,1]}s_{\al}\for a=\al^{\vee},
\ \al\in R.
&\eqnu
}
$$

The {\it extended Weyl group} $ W^b$ generated by $W\and B'$
(instead of $A'$) is isomorphic to $W\lsmash B'$:
$$
\eqalignno{
&(wb')([z,\ze])\ =\ [w(z),\ze-(z,b)] \for w\in W, b\in B.
&\eqnu
}
$$

 Given $b_+\in B_+$, let
$$
\eqalignno{
&\om_{b_+} = w_0w^+_0  \in  W,\ \pi_{b_+} =
b'_+(\om_{b_+})^{-1}
\ \in \ W^b, \ \om_i=\om_{b_i},\pi_i=\pi_{b_i},
&\eqnu
\label\w0\eqnum*
}
$$
where $w_0$ (respectively, $w^+_0$) is the longest element in $W$
(respectively, in $ W_{b_+}$ generated by $s_i$ preserving $b_+$)
relative to the
set of generators $\{s_i\}$ for $i >0$.

We will  use here only the
elements $\pi_r=\pi_{b_r}, r \in O$. They leave $\Ga^a$ invariant
and form a group denoted by $\Pi$,
 which is isomorphic to $B/A$ by the natural
projection $\{b_r \to \pi_r\}$. As to $\{\om_r\}$,
they preserve the set $\{-\th,\al_i, i>0\}$.
The relations $\pi_r(\al_0)= \al_r= (\om_r)^{-1}(-\th)
$ distinguish the
indices $r \in O^*$. Moreover (see e.g. [C2]):
$$
\eqalignno{
& W^b  = \Pi \lsmash W^a, \where
  \pi_rs_i\pi_r^{-1}  =  s_j \if \pi_r(\al_i)=\al_j,\  0\le j\le n.
&\eqnu
}
$$

We extend the notion
of the length to $W^b$.
Given $\nu\in\nu_R,\  r\in O^*,\  \tw \in W^a$, and a reduced
decomposition $\tw\ =\ s_{j_l}...s_{j_2} s_{j_1} $ with respect to
$\{s_j, 0\le j\le n\}$, we call $l\ =\ l(\hw)$ the {\it length} of
$\hw = \pi_r\tw \in W^b$. Setting
$$
\eqalign{
\la(\hw) = &\{ \tal^1=\al_{j_1},\
\tal^2=s_{j_1}(\al_{j_2}),\
\tal^3=s_{j_1}s_{j_2}(\al_{j_3}),\ldots \cr
&\ldots,\tal^l=\tw^{-1}s_{j_l}(\al_{j_l}) \},
}
\eqnu
$$
\label\tal\eqnum*
one can represent
$$
\eqalign
{
&l=|\la(\hw)|=\sum_\nu l_\nu, \for l_\nu = l_\nu(\hw)=|\la_\nu(\hw)|,\cr
&\la_\nu(\hw) = \{\tal^{m},\ \nu(\tal^{m})= \nu(\tal_{j_m})= \nu\},
1\le m\le l,
}
\eqnu
\label\laset\eqnum*
$$
where $|\ |$  denotes the  number of elements,
 $\nu([\al,k]) \equal \nu_{\al}$.

For instance,
$$
\eqalign{
&l_\nu(b')\ =\ \sum_{\al} |(b,\al)|,\  \al\in R_+,
\nu_\al=\nu \in \nu_R, \cr
&l_\nu(b_+')\ =\ 2(b_+,\rho_\nu) \when b_+ \in B_+.
}
\eqnu
$$
\label\lb\eqnum*
Here $|\ | = $
absolute value. Later on  $b$ and $b'$ will not be distinguished.

We put
$m=2 \for D_{2k} \and C_{2k+1},\ m=1 \for C_{2k}, B_{k}$,
otherwise $m=|\Pi|$. The definition involves
the parameters $ q,\{ t_\nu , \nu \in \nu_R \}$ and independent
variables $X_1,\ldots,X_n$.
Let us set
$$
\eqalignno{
&   t_{\tal} = t_{\nu(\tal)},\ t_j = t_{\al_j},
\where \tal \in R^a, 0\le j\le n, \cr
& X_{\tb}\ =\ \prod_{i=1}^nX_i^{k_i} q^{ k}
\if \tb=[b,k],
&\eqnu \cr
&\for b=\sum_{i=1}^nk_i b_i\in B,\ k \in {1\over m}\Z.
}
$$
\label\Xde\eqnum*

 Later on $ \C_{ q}$
 is the field of rational
functions in $ q^{1/m},$
$\C_q[X] = \C_q[X_b]$  means the algebra of
polynomials in terms of $X_i^{\pm 1}$
with the coefficients depending
on $ q^{1/m}$ rationally. We replace $ \C_{ q}$
by $ \C_{ q,t}$ if the functions (coefficients)
also depend rationally
on $\{t_\nu^{1/2} \}$.

Let $([a,k],[b,l])=(a,b)$ for $a,b\in B,\
 a_0=\al_0,\ \nu_{\al^\vee}=\nu_\al, $
and  $ \al_{r^*} \equal \pi_r^{-1}(\al_0)$ for $  r\in O^*$.

\proclaim{Definition }
 The  double  affine Hecke algebra $\HH\ $
(see [C1,C2])
is generated over the field $ \C_{ q,t}$ by
the elements $\{ T_j,\ 0\le j\le n\}$,
pairwise commutative $\{X_b, \ b\in B\}$ satisfying (\ref\Xde),
 and the group $\Pi$ where the following relations are imposed:

(o)\ \  $ (T_j-t_j^{1/2})(T_j+t_j^{-1/2})\ =\ 0,\ 0\ \le\ j\ \le\ n$;

(i)\ \ \ $ T_iT_jT_i...\ =\ T_jT_iT_j...,\ m_{ij}$ factors on each side;

(ii)\ \   $ \pi_rT_i\pi_r^{-1}\ =\ T_j \if \pi_r(\al_i)=\al_j$;

(iii)\  $T_iX_b T_i\ =\ X_b X_{a_i}^{-1} \if (b,\al_i)=1,\
1 \le i\le  n$;

(iv)\  $T_0X_b T_0\ =\ X_{s_0(b)}\ =\ X_b X_{\th} q^{-1}
\if (b,\th)=-1$;

(v)\ \ $T_iX_b\ =\ X_b T_i$ if $(b,\al_i)=0 \for 0 \le i\le  n$;

(vi)\ $\pi_rX_b \pi_r^{-1}\ =\ X_{\pi_r(b)}\ =\ X_{\om^{-1}_r(b)}
 q^{(b_{r^*},b)},\  r\in O^*$.
\label\double\theoremnum*
\endproclaim
\proofbox

Given $\tw \in W^a, r\in O,\ $ the product
$$
\eqalignno{
&T_{\pi_r\tw}\equal \pi_r\prod_{k=1}^l T_{i_k},\where
\tw=\prod_{k=1}^l s_{i_k},
l=l(\tw),
&\eqnu
\label\Tw\eqnum*
}
$$
does not depend on the choice of the reduced decomposition
(because $\{T\}$ satisfy the same ``braid'' relations as $\{s\}$ do).
Moreover,
$$
\eqalignno{
&T_{\hv}T_{\hw}\ =\ T_{\hv\hw}\  \hbox{ whenever}\
 l(\hv\hw)=l(\hv)+l(\hw) \for
\hv,\hw \in W^b.
&\eqnu}
$$
\label\TT\eqnum*
  In particular, we arrive at the pairwise
commutative elements
$$
\eqalignno{
& Y_{b}\ =\  \prod_{i=1}^nY_i^{k_i} \if
b=\sum_{i=1}^nk_ib_i\in B,\where
 Y_i\equal T_{b_i},
&\eqnu
\label\Yb\eqnum*
}
$$
satisfying the relations
$$
\eqalign{
&T^{-1}_iY_b T^{-1}_i\ =\ Y_b Y_{a_i}^{-1} \if (b,\al_i)=1,
\cr
& T_iY_b\ =\ Y_b T_i \if (b,\al_i)=0, \ 1 \le i\le  n.}
\eqnu
$$
Let us introduce the following elements from
$\C_t^n$:
$$
\eqalign{
&t^{\pm\rho}\equal (l_t(b_1)^{\pm 1},\ldots,l_t(b_n)^{\pm 1}),\where\cr
&l_t(\hw)\equal \ \prod_{\nu\in\nu_R} t_\nu^{l_\nu(\hw)/2},\
\hw\in W^b,
}
\eqnu
\label\qlen\eqnum*
$$
and the corresponding {\it evaluation maps}:
$$
\eqalign{
&X_i(t^{\pm\rho})= l_t(b_i)^{\pm 1} = Y_i(t^{\pm\rho}),\ 1\le i\le n.
}
\eqnu
\label\eval\eqnum*
$$
For instance, $X_{a_i}(t^{\rho})\ =\ l_t(a_i)= t_i$ (see (\ref\lb)).

\proclaim {Theorem}
i) The elements $H \in \HH\ $  have
the unique decompositions
$$
\eqalignno{
&H =\sum_{w\in W }  g_{w}  T_{w} f_w,\
g_{w} \in \C_{ q,t}[X],\ f_{w} \in \C_{ q,t}[Y].
&\eqnu
}
$$

ii) The   map
$$
\eqalign{
 \vph: &X_i \to Y_i^{-1},\ \  Y_i \to X_i^{-1},\  \ T_i \to T_i, \cr
&t_\nu \to t_\nu,\
 q\to  q,\ \nu\in \nu_R,\ 1\le i\le n.
}
\eqnu
\label\vph\eqnum*
$$
can be extended to an anti-involution
($\vph(AB)=\vph(B)\vph(A)$)
 of \HH\ .

iii) The  linear functional on \HH\
$$
\eqalignno{
&[\![ \sum_{w\in W }  g_{w}  T_{w} f_{w}]\!]\ =\
\sum_{w\in W} g_{w}(t^{-\rho}) l_t(w) f_{w}(t^{\rho})
&\eqnu
\label\brack\eqnum*
}
$$
is invariant with respect to $\vph$. The bilinear form
$$
\eqalignno{
&[\![ G,H]\!]\equal [\![ \vph(G)H]\!],\
G,H\in \HH\ ,
&\eqnu
\label\form\eqnum*
}
$$
is symmetric ($[\![ G,H]\!]= [\![ H,G]\!]$)
and non-degenerate.
\endproclaim
\label\dual\theoremnum*

{\it Proof.} The first statement is  from Theorem 2.3 [C2].
The map $\vph$ is the composition of the
involution (see [C1])
$$
\eqalign{
  \vep:\ &X_i \to Y_i,\ \  Y_i \to X_i,\  \ T_i \to T_i^{-1}, \cr
&t_\nu \to t_\nu^{-1},\
 q\to  q^{-1},\ 1\le i\le n,
}
\eqnu
\label\vep\eqnum*
$$
and the main anti-involution  from [C2], introduced by
the condition
$\langle Hf,g\rangle  = $ $ \langle f,H^*g\rangle$
 (see (\ref\innerpro)):
$$
\eqalign{
  & X_i^*\ =\  X_i^{-1},\   Y_i^*\ =\  Y_i^{-1},\
 T_i^* \ =\  T_i^{-1}, \cr
&t_\nu \to t_\nu^{-1},\
 q\to  q^{-1},\ 0\le i\le n.
}
\eqnu
\label\star\eqnum*
$$
The other claims follow directly from the definition of
$[\![ \ ]\!]$.
\proofbox

Let us give the explicit formulas  for the action of $\vph,\vep$
on $T_0$:
$$
\eqalign{
&\vph(T_0)\ =\ Y_\th^{-1}T_0X_\th^{-1}\ =\ T_{s_\th}^{-1}X_\th^{-1},\cr
&\vep(T_0)\ =\ X_\th T_0^{-1}Y_\th\ =\ X_\th T_{s_\th}.
}
\eqnu
\label\vpht0\eqnum*
$$

One can extend $[\![\ ]\!]$ to the localization of \HH\
with respect to all polynomials in $X$ (or in $Y$). The
algebra becomes the semi-direct product of $\C[W^b]$ and
$\C(X)$ after this (see [C3]). Sometimes it is also convenient
to involve proper completions of $\C(X)$.

%
%
%
%
\section { Difference operators}
Setting (see  Introduction)
$$
\eqalignno{
& x_{\tb}=  \prod_{i=1}^nx_i^{k_i} q^{ k} \if
\tb=[b,k],
b=\sum_{i=1}^nk_i b_i\in B,\ k \in {1\over m}\Z,
&\eqnu
\label\xde\eqnum*}
$$
for independent $x_1,\ldots,x_n$, we will
 consider $\{X\}$ as  operators acting in $\C_q[x]=
\C_q[x_1^{\pm 1},\ldots,x_n^{\pm 1}]$:
$$
\eqalignno{
& X_{\tb}(p(x))\ =\ x_{\tb} p(x),\    p(x) \in
\C_q [x].
&\eqnu}
$$
\label\X\eqnum*
The elements $\hw \in W^b$ act in $\C_{ q}[x]$
 by the
formulas:
$$
\eqalignno{
&\hw(x_{\tb})\ =\ x_{\hw(\tb)}.
&\eqnu}
$$
 In particular:
$$
\eqalignno{
&\pi_r(x_{b})\ =\ x_{\om^{-1}_r(b)} q^{(b_{r^*},b)}
\for \al_{r^*}\ =\ \pi_r^{-1}(\al_0), \ r\in O^*.
&\eqnu}
$$
\label\pi\eqnum*

The {\it Demazure-Lusztig operators} (see
[KL1, KK, C1], and [C2] for more detail )
$$
\eqalignno{
&\hT_j\  = \  t_j ^{1/2} s_j\ +\
(t_j^{1/2}-t_j^{-1/2})(X_{a_j}-1)^{-1}(s_j-1),
\ 0\le j\le n.
&\eqnu
\label\Demaz\eqnum*
}
$$
act   in $\C_{ q,t}[x]$ naturally.
We note that only $\hT_0$ depends on $ q$:
$$
\eqalign{
&\hT_0\  =  t_0^{1/2}s_0\ +\ (t_0^{1/2}-t_0^{-1/2})
( q X_{\th}^{-1} -1)^{-1}(s_0-1),\cr
&\where
s_0(X_i)\ =\ X_iX_{\th}^{-(b_i,\th)} q^{(b_i,\th)}.
}
\eqnu
$$

\proclaim{Theorem }
 The map $ T_j\to \hT_j,\ X_b \to X_b$ (see (\ref\Xde,\ref\X)),
$\pi_r\to \pi_r$  (see (\ref\pi)) induces a $ \C_{ q,t}$-linear
homomorphism from \HH\ to the algebra of linear endomorphisms
of $\C_{ q,t}[x]$.
 This representation is faithful and
remains faithful when   $  q,t$ take  any non-zero
values assuming that
 $ q$ is not a root of unity (see [C2]). The image $\hat{H}$
is uniquely determined from the following condition:
$$
\eqalign{
&\hat{H}(f(x))\ =\ g(x)\for H\in \HH\ ,\if Hf(X)-g(X)\ \in\cr
 &I\equal \{\sum_{i=0}^n H_i(T_i-t_i)+
\sum_{r\in O^*} H_r(\pi_r-1), \where H_i,H_r\in \HH\ \}.
}
\eqnu
\label\hat\eqnum*
$$
\endproclaim
\proofbox
\label\faith\theoremnum*

Due to  Theorem \ref\dual, an arbitrary  $H\in \HH\ $ can be
uniquely represented in the form
$$
\eqalign{
H =&\sum_{b\in B, w\in W }  g_{b,w}  Y_b T_{w},\
g_{b,w} \in \C_{ q,t}[X],\cr
=&\sum_{b\in B, w\in W }   T_{w} X_b  g'_{b,w},\
g'_{b,w} \in \C_{ q,t}[Y].
}
\eqnu
$$
We set:
$$
\eqalign{
&[H]_{\dagger}\  =\ \sum_{b\in B, w\in W }  g_{b,w}  Y_b l_t(w),\
{}_{\dagger}[H]\ =\ \sum_{b\in B, w\in W }   l_t(w) X_b  g'_{b,w},\cr
&[H]_\ddagger\  =\ \sum_{b\in B, w\in W }  g_{b,w}  [\![Y_bT_w]\!],\
{}_\ddagger[H]\ =\ \sum_{b\in B, w\in W }    [\![T_wX_b]\!]  g'_{b,w}.
}
\eqnu
\label\redH\eqnum*
$$
One easily checks that
$$
\eqalign{
[\![H_1 H_2]\!]\ =\ &[\![H_1 [H_2]_{\dagger}]\!]\ =\
[\![ {}_{\dagger}[H_1]H_2]\!]\ =\cr
&[\![H_1 [H_2]_\ddagger]\!]\ =\ [\![ {}_\ddagger[H_1]H_2]\!]
\for H_1,H_2\in \HH\ .
}
\eqnu
\label\dag\eqnum*
$$

The image $\hH$ of $H$ can be uniquely
represented as follows:
$$
\eqalignno{
&\hH\ = \sum_{b\in B, w\in W} h_{b,w} b  w
\ =\ \sum_{b\in B, w\in W} w b h'_{b,w}.
&\eqnu
\label\hatH\eqnum*
}
$$
where  $h_{b,w}, h'_{b,w}$ belong to  the field $\C_{ q,t}(X)$
of rational  functions in
$X_1,...,X_n$. We extend the above operations to arbitrary operators
in the form (\ref\hatH):
$$
\eqalign{
&[\hH]_{\dagger}= \sum h_{b,w}  b ,\
{}_{\dagger}[\hH] = \sum b h'_{b,w}\  ,\
[\![ \hH ]\!]\
 = \sum h_{b,w}
(t^{-\rho}).
}
\eqnu
\label\Brack\eqnum*
$$
These operations commute with the homomorphism $H\to\hH$.

Let us define the
{\it difference Harish-Chandra map} (see
[C2], Proposition 3.1):
$$
\eqalignno{
&\chi( \sum_{w\in W,b\in B}  h_{b,w} b w)\ =\
\sum_{b\in B,w\in B} h_{b,w}(\diamondsuit)
y_b \ \in \C_{ q,t}[y],
&\eqnu
}
$$
where
$\diamondsuit \equal ( X_1=...=X_n=0 ),\ \{y_b\}$ is one more  set of variables
 introduced  for
independent $y_1,...,y_n$ in the same way as
$\{x_b\}$ were.

\proclaim{Proposition}
 Setting
$$
\eqalign{
&\l_f \ =\ f(Y), \ \hat{\l}_f=\ f(\hY),\
L_f = L_f^{ q,t}\ \equal\ [(\hat{\l}_f)]_{\dagger}
}
\eqnu
\label\Ll\eqnum*
$$
for $f=\sum_b g_b y_b \in \C_{ q,t}[y]$, one has:
$$
\eqalign{
\chi(\hat{\l}_f)\ =\  \chi(L_f)\ =\
[\![ f(Y)]\!]\ =\
\sum_{b\in B}  g_{b}
\prod_\nu t_\nu^{(b,\rho_\nu)} y_b.
}
\eqnu
$$
\label\chi\eqnum*
\endproclaim
\proofbox

The proof
of the following theorem repeats the proof of Theorem 4.5,[C2]
(where the relations $t_\nu= q_\nu^{k_\nu}$ for $k_\nu\in \Z_+$
were imposed). We note that once  (\ref\Lf) is known for these
special $t$ it holds true for all $ q,t$ since all the
coefficients of difference operators and polynomials
are rational in $ q,t$.

\proclaim{Theorem}
The difference operators $\{ L_f, \
f(y_1,\ldots,y_n)\in \C_{ q,t}[y]^W\}$
are pairwise commutative,  $W$-invariant (i.e $w L_f w^{-1}=$
$L_f$ for all $w\in W$) and preserve $\C_{ q,t}[x]^W$. The
Macdonald polynomials $p_b=p_b^{ q,t} (b\in B_-)$
from  (\ref\macd) are their eigenvectors:
$$
\eqalignno{
&L_f(p_b^{ q,t})=f(t^\rho q^{-b}) p_b^{ q,t},\
y_i(t^\rho q^{-b})\equal
 q^{-(b_i,b)}\prod_\nu t_\nu^{(b_i,\rho_\nu)}.
&\eqnu
\label\Lf\eqnum*}
$$
\label\LF\theoremnum*
\endproclaim
\proofbox

We fix a subset $v\in \nu_R$ and introduce the
{\it shift operator} by the formula
$$
\eqalignno{
&\g_v \ =\
(\x_v) ^{-1}\y_v,\ G_v^{ q,t}\ =\ [\hat{\g}_v]_{\dagger}  \ =
\ (\x_v )^{-1}[\y_v]_{\dagger},
 &\eqnu
\label\shift\eqnum*
}
$$
$$
 \x_v  = \prod_{\nu_a\in v}((t_a X_{a})^{1/2}-
(t_a X_{a})^{-1/2}),\  \y_v  = \prod_{\nu_a\in v}
(t_a Y_{a}^{-1})^{1/2}-
(t_a Y_{a}^{-1})^{-1/2}).
$$
Here $a=\al^\vee\in R_+^\vee, \nu_{a}=\nu_\al, t_a=t_\al$, the
elements $\x_{v}= \x_{v}^t, \y_{v}=\y_{v}^t $
belong to $\C_t [X],
\C_t [Y]$ respectively.

\proclaim{Theorem}
The operators $\hat{\g}_{v}$
and $G_v^{ q,t}$  are $W$-inva\-riant and preserve $ \C_{ q,t}[x]^W$
(their restrictions to the latter space coincide). Moreover,
if $t_{\nu}=1$ when
$\nu\not\in v$ then
$$
\eqalign{
& G_v^{ q,t}L_f^{ q,t}\ =\ L_f^{ q, tq_v} G_v^{ q,t}
\for f\in \C_{ q,t}[y]^W,
\cr
&G_v^{ q,t} (p_{b}^{ q,t})= g_v^{ q,t}(b)
p_{b+r_v}^{ q, tq_v}, \for\cr
&g_v^{ q,t}(b)\ =\
\prod_{a\in R_+^\vee,\nu_a\in v} (y_a(t^{\rho/2} q^{-b/2}) -
t_a y_a(t^{-\rho/2} q^{+b/2})),
}
\eqnu
\label\Gnu\eqnum*
$$
where $r_v=\sum_{\nu\in v}r_\nu,\ tq_v=\{t_\nu q^{2/\nu} , t_{\nu'}\}$
for $\nu\in v\not\ni \nu'\ $, $p_c=0 \for c\not\in B_-$.
\label\Gp\theoremnum*
\endproclaim
{\it Proof.}  When $t_\nu= q^{2k_\nu/\nu}$ for $k_\nu\in \Z_+$
these statements  are in fact from [C2].
They give (\ref\Gnu) for all $ q,t$. Indeed,
 it can be rewritten as follows:
$$
\eqalign{
& [\hat{\l}_f^{ q,t}\x_v^{t}]_\dagger\  =\
\x_v^{t} L_f^{ q,tq_v},
}
\eqnu
\label\xGx\eqnum*
$$
where the coefficients of the difference operators on both
sides are from
$\C_{ q,t}[X]$.
Here we used that $[\l\m]_{\dagger}=[\l]_{\dagger}[\m]_{\dagger} $ for
arbitrary operators  $\l,\m$ in the form (\ref\hatH) if the
second is $W$-invariant. The remaining formulas can be
deduced from [C2] in the same way (they mean certain
identities in $\C_{ q,t}$ which are enough to check for
$t_\nu= q^{2k_\nu/\nu}$). One can use (\ref\chi) as well.
\proofbox

We will also need Proposition 3.4 from [C2]:
\proclaim{Proposition}
Given $b\in B_-,$ let $\ m_b=  \sum_{w\in W/W_b} x_{w(b)}$
for the stabilizer $W_b$ of $b$ in $W$,
$L_b\equal L_{m_b}  =  \{L_{b}\} + \sum_{c } f_{b}^c(X)c$,\
$ W(c) \succ b$. Then
$$
\eqalign{
& \{L_{b}\} = \sum_{w\in W/W_b}
 \prod _{a\in R_+^\vee, j}
{t^{1/2}_{a}X_{w(a)} q_a^j - t_{a}^{-1/2}\over
X_{w(a)} q_a^j- 1} w(b),\
-(b,a^\vee)>j\ge 0.
}
\eqnu
\label\lead\eqnum*
$$
 If $r\in O^*$ then $L_{-b_r}\  = \ \{L_{-b_r}\}$.
\label\LEAD\theoremnum*
\endproclaim

%
%
%
%
\section { Duality and evaluation conjectures}
First of all we will use Theorem \ref\dual to define the
{\it
Fourier pairing}. In the classical theory the latter is the
 inner product of a function and the
Fourier transform  of another function.
In this and the next sections we will
identify the elements
$H\in \HH\ $ with their images $\hH$. The following pairing
on $f,g\in \C_{ q,t}[x]$ is symmetric and non-degenerate:
$$
\eqalign{
 &[\![f,g]\!]= [\![f(X),g(X)]\!] = [\![\vph(f(X))g(X)]\!]  =\cr
&[\![\bar{f}(Y)g(X)]\!] =  \{\l_{\bar{f}}(g(x))\}(t^{-\rho}),\cr
&\bar{x}_b\ =\ x_{-b}\ =\ x_b^{-1},\ \bar{ q}\ =\  q,\
\bar{t}\ =\ t,
}
\eqnu
\label\Fourier\eqnum*
$$
where
${\l}$ is from (\ref\Ll),
and we used the main defining property (\ref\hat)
of the representation
from  Theorem \ref\faith. The pairing remains non-degenerate
when restricted to $W$-invariant polynomials.

\proclaim{ Definition}
The Fourier adjoints $\vph(\l),\vph(L)$ of $\C_{ q,t}$-linear
operators $\l,L$
acting respectively  in $\C_{ q,t}[x]$ or
in $\C_{ q,t}[x]^W$ are defined from the relations:
$$
\eqalign{
&[\![\l(f),g]\!]\ =\  [\![f,\vph(\l)(g)]\!],\ f,g\in C_{ q,t}[x],\cr
&[\![L(f),g]\!]\  =\   [\![f,\vph(L)(g)]\!],\ f,g\in \C_{ q,t}[x]^W.
}
\eqnu
$$
If $\l$ preserves $\C_{ q,t}[x]^W$ then so does $\vph(\l)$
and $\vph(L)=
[\vph(\l)]_{\dagger}, \where L=[\l]_{\dagger}$ is the restriction of
$\l$ to the invariant polynomials.
\endproclaim
\label\FT\theoremnum*
\proofbox

This anti-involution  ($\vph^2=\hbox{id}$) extends  $\vph$ from
(\ref\vph) by construction. If $f\in C_{ q,t}[x]^W$, then
$\vph(L_f)= [\bar{f}(X)]_{\dagger}$. We arrive at the following
theorem:

\proclaim {Duality Theorem}
Given $b,c\in B_-$ and the corresponding Macdonald's
polynomials $p_b, p_c$,
$$
p_b(t^{-\rho}  q^{c})p_c(t^{-\rho})\ =
[\![p_b,p_c]\!]\ =\ [\![p_c,p_b]\!]\ = \
p_c(t^{-\rho}  q^{b})p_b(t^{-\rho}).
\eqnu
\label\pp\eqnum*
$$
\endproclaim
\proofbox

To complete this theme we need to calculate $p_b(t^{-\rho})$.
The main step is the formula for
$p'((tq_v)^{-\rho})$ in terms of $p(t^{-\rho})$,
where  (see (\ref\Gnu))
$$p=p_b,\ p'\ =\ p_{b+r_v}^{tq_v},\
p'=(g_v^{t}(b))^{-1}
G_v^{t} (p).$$
Here and in similar formulas we show the dependence on $t$
omitting $ q$
since the latter  will be the same for all polynomials and
operators. Let
$$ \bar{\y}_v^t  = \prod_{a\in R_+^\vee,\nu_a\in v}
((t_a Y_{a})^{1/2}-
(t_a Y_{a})^{-1/2}).$$

\proclaim { Key Lemma}
$$
\eqalignno{
d_v^t p'((tq_v)^{-\rho})& =
\prod_{a\in R_+^\vee,\nu_a\in v} \Bigl(t_a^{-1} y_a(t^{-\rho/2} q^{+b/2})
-y_a(t^{+\rho/2} q^{-b/2})\Bigr)p(t^{-\rho}), \cr
d_v^t& =
\prod_{a\in R_+^\vee,\nu_a\in v} \Bigl(t_a^{-1} y_a((tq_v)^
{-\rho/2})
-y_a((tq_v)^{+\rho/2})\Bigr)m_{-r_v}(t^{-\rho}).&\eqnu
\label\key\eqnum*
}
$$
\label\KEY\theoremnum*
\endproclaim
{\it Proof.}
Let us use formula (\ref\xGx):
$$
\eqalign{
&[\![({\y}_v^t\x_v^t)(\x_v^t)^{-1}\l^t_{\bar{p}'}\x_v^t]\!]\ =
\ [\![({\y}_v^t\x_v^t)[(\x_v^t)^{-1}\l^t_{\bar{p}'}\x_v^t]_{\dagger}
]\!]\ =\cr
&[\![(\hat{\y}_v^t\hat{\x}_v^t)   \bar{p}'(\hat{Y}^{tq_v})]\!]\ =
\ [\![({\y}_v^t\x_v^t) ]\!] p'((tq_v)^{-\rho}).
}
\eqnu
\label\one\eqnum*
$$
On the other hand, it equals:
$$
\eqalign{
&[\![{\y}_v^t p'(Y) \x_v^t ]\!]\ =
[\![\y_v^t {p}'({X}^{t}) {\x}_v^t ]\!]\ =\cr
&[\![{\y}_v^t (\x_v^t p'(x)) ]\!]\ =
\pm[\![\bar{\y}_v^t (\x_v^t p'(x)) ]\!].
}
\eqnu
$$
Here we applied the anti-involution $\vph\  (\vph(\x)=\y,\
\vph(\y)=\x)$, then went
 from the abstract
$[\![\ ]\!]$ to that from (\ref\Brack), and
used  Theorem \ref\faith. The last transformation requires
 special comment. We will justify it in a moment.

After this, one can use  (\ref\Lf):
$$
\eqalign{
&[\![\bar{\y}_v^t (\x_v^t p'(x)) ]\!]\ =
[\![(\bar{\y}_v^t {\y}_v^t)   g_v^t(b)^{-1}p(x)) ]\!]\ =\cr
& g_v^t(b)^{-1}(\bar{\y}\y)(t^\rho q^{-b})
[\![ p(x)]\!]\ =\cr
&\prod_{a\in R_+^\vee,\nu_a\in v} (t_a^{-1} y_a(t^{-\rho/2} q^{+b/2}) -
y_a(t^{+\rho/2} q^{-b/2}))[\![ p(x) ]\!].
}
\eqnu
\label\evalp\eqnum*
$$
Finally, ${d}_v^t\equal \pm[\![{\y}_v^t\x_v^t ]\!]$ can be
determined from (\ref\evalp) and the relation
$1=p'\ =\ g_v^{t}(b)^{-1}
G_v^{t} (p_{b}^{t})\for b=-r_v,$ where $p_{-r_v}$ coincides with
the monomial function $m_{-r_v}$ (it follows directly from
the definition):
$$
\eqalign{
&{d}_v^t\ =\
\prod_{a\in R_+^\vee,\nu_a\in v} (t_a^{-1} y_a((tq_v)^
{-\rho/2})
-y_a((tq_v)^{+\rho/2})) m_{-r_v}(t^{-\rho}).
}
\eqnu
\label\dterm\eqnum*
$$

Let us check that
$$[\![(\bar{\y}_v^t-l_\ep(w_0){\y}_v^t)
(\x_v^t p'(x)) ]\!]=0 \for \hbox {\ any\ }
p'\in \C[x],
$$
where $l_\ep(w_0)= \prod_\nu \ep_\nu^{l_\nu(w_0)}$,
$$
\ep =   \{\ep_\nu= -1 \if \nu\in v,\hbox{\ otherwise\ } \ep_\nu=1\},\
 \nu\in \nu_R.
$$
Following formula (4.18),[C2] we introduce
the
{\it $t$-symmetrizers}, setting
$$
\eqalign{
&\p_v^t\ =\ (\pi_v^t)^{-1}\sum_{w\in W}
\prod_\nu(\ep_\nu t_\nu^{1/2})^{\ep_\nu(l_\nu(w)-l_\nu(w_0))} T_w,
\cr
&\pi_v^t\ =\ \sum_{w\in W}
\prod_\nu (\ep_\nu t_\nu^{1/2})^{\ep_\nu(2l_\nu(w)-l_\nu(w_0))}.
}
\eqnu
$$
It results from Proposition 3.5 and  Corollary 4.7(ibidem) that
$$\p_v^t (\x_v^t p')=\x_v^t p',\  \hat{\p}_v^t\p_v^{t=0}=\hat{\p}_v^t, \if
 t_\nu=1 \for
\nu\not\in \nu_R. $$
Hence
$$
\eqalign{
& [\![(\bar{\y}_v^t -l_\ep(w_0){\y}_v^t) (\x_v^t p'(x)) ]\!] =
[\![(\bar{\y}_v^t -l_\ep(w_0){\y}_v^t)\p_v^t  (\x_v^t p'(x)) ]\!] =\cr
&[\![ (\y_v^t \bar{p}'(Y))\p_v^t (\bar{\x}_v^t-l_\ep(w_0){\x}_v^t) ]\!] =
[\![ (\y_v^t \bar{p}'(Y))
\{\hat{\p}_v^t\p_v^{t=0} (\bar{\x}_v^t-l_\ep(w_0){\x}_v^t)\}
 ]\!].
}
$$
The latter equals zero.
\proofbox

Let us take any set $k=\{k_{\nu_1}\ge k_{\nu_2}\}\in \Z_+$ and put
$$
\eqalign{
t(k)=\{ q^{2k_\nu/\nu}\},\   k\cdot r=\sum_\nu k_\nu r_\nu,\
p_b^{(k)}= p_b^{t(k)}.
}
\eqnu
$$
The remaining  part of the calculation is based on
 the following chain of the shift operators that will be applied
to $p^{(0)}_{b-k\cdot r}=m_{b-k\cdot r}$
one after another:
$$
\eqalignno{
&G_{\nu_R}^{(k-1)}G_{\nu_R}^{(k-2)}\cdots
G_{\nu_R}^{(k-s)}G_{\nu_1}^{(k-s-e)}
\cdots G_{\nu_1}^{(0)},
&\eqnu
}
$$
where
$ k_{\nu_1}=s+d,\ k_{\nu_2}=s,\
e=\{e_\nu\},\  e_{\nu_1}=1,\ e_{\nu_2}=0,\ k-s= de$,
the set $\{1,1\}$ is denoted by $1$.

Lemma \ref\KEY gives that for a certain $D^{(k)}$ (which does
not depend on $b$):
$$
\eqalignno{
&D^{(k)}p_b^{(k)}(t(k)^{-\rho})\ =  &\eqnu \cr
m_{b-k\cdot r}(1)&\prod_{a\in R_+^\vee,\nu_a\in v(i)}^{0\le i<s+d}
\Bigl(t(i)_a\, y_a(t(i)^{\rho/2} q^{-b(i)/2})
-y_a(t(i)^{-\rho/2} q^{b(i)/2})\Bigr),
\label\prodb\eqnum*
}
$$
where $t(i)=t(k(i)),\ b(i)=b -(k-k(i))\cdot r$,
$$k(i)=ie,\ v(i)=\nu_1\if i < d, \
 k(i)= i-d + de, \
v(i)=\nu_R \if i\ge d.$$

As to $D^{(k)}$, it equals the right hand side of (\ref\prodb)
when $b=0$. We note that
$t(i)_a=  q^{2j/\nu_a} \for j=k_a+i-s-d,\ k_a=k_{\nu_a},$
because  $i\ge d\if \nu_a\neq \nu_1$ (and $0\le j< k_a$).
The relation  $(2/\nu)\rho_\nu=r_\nu$ leads to the formulas:
$$
\eqalign{
&t(i)^{\rho/2} q^{-b(i)/2}\ =\
 q^{(k(i)\cdot r -b+(k-k(i))\cdot r)/2}\ =\   q^{k\cdot r-b/2},\cr
&y_a(t(i)^{\rho/2} q^{-b(i)/2})\ =\
 q^{ (k\cdot r-b,a)/2 }.
}
\eqnu
$$
 Finally, we arrive at the following theorem:
\proclaim{ Evaluation Theorem}
$$
\eqalignno{
&p_b^{(k)}(t(k)^{-\rho})= \ &\eqnu\cr
{m_{b-k\cdot r}(1)\over m_{-k\cdot r}(1)}
&\prod_{\al\in R_+, 0\le j< k_\al}
\Bigl(
{
 q^{ \{(k\cdot r-b,\al)+j\}/\nu_\al }
- q^{- \{(k\cdot r-b,\al)+j\}/\nu_\al}
 \over
 q^{ \{(k\cdot r,\al)+j\}/\nu_\al }
- q^{- \{(k\cdot r,\al)+j\}/\nu_\al}
}
\Bigr).
\label\evde\eqnum*
}
$$
\endproclaim
\proofbox

Here
$m_{b-k\cdot r}(1)/ m_{-k\cdot r}(1)\ =\
|W(b-k\cdot r)|/|W(k\cdot r)|\in \Z_+$. It  equals $1$
for all $b\in B_-$ when
$\prod_\nu k_\nu\neq 0$. Assuming this
we have:
$$
\eqalignno{
&p_b^{t(k)}(t(k)^{-\rho})\ =  &\eqnu \cr
 q^{(k\cdot r,b)}&
\prod_{\al\in R_+, 0\le j< \infty}
\Bigl(
{
(1- q^{ 2\{(k\cdot r-b,\al)+j\}/\nu_\al })
(1-t_\al q^{ 2\{(k\cdot r,\al)+j\}/\nu_\al })
 \over
(1-t_\al q^{ 2\{(k\cdot r-b,\al) +j\}/\nu_\al })
(1- q^{ 2\{(k\cdot r,\al)+j\}/\nu_\al })
}
\Bigr).
\label\evdeq\eqnum*
}
$$
The limit of (\ref\evdeq) as one of the $k_\nu$ approaches zero
exists and coincides with (\ref\evde). Since both sides of
this formula are rational functions in $t(k)\and  q$ we get
 (\ref\EV)  (cf.\ Theorem \ref\Gp).

We note that actually this paper does not depend very much
on the definition of the Macdonald polynomials from the
Introduction. We can eliminate $\mu$ introducing these
polynomials as the eigenfunctions of the $L$-operators
(formula (\ref\Lf)). Therefore it is likely that
paper [C4] can be extended to give  a "difference-elliptic"
Weyl dimension formula.
%
%
%
%
\section { Discretization, applications}
Continuing the same line let us establish the {\it
recurrence relations} for the Macdonald polynomials
generalizing the three-term relation for the $q$-ultraspherical
polynomials (Askey, Ismail) and the Pieri rules.
We need to go to the lattice
version of the considered functions and operators. The {\it
 discretization} of
functions $g(x),\ x\in \C^n$  and  ($ q$-)difference operators
is defined
as follows:
$$
\eqalign{
& {}^\de g(b)\ =\ g( q^b t^{-\rho}),\
({}^\de a ({}^\de g))(b)\ =\ {}^\de g(b-a),\  a,b\in B, \cr
&({}^\de X_a({}^\de g))(b)\ =\  x_{a}( q^{b}t^{-\rho})
\ {}^\de g(b)\ = \
 q^{(a,b)}\prod_\nu t_\nu^{-(a,\rho_\nu)}\ {}^\de g(b).
}
\eqnu
\label\deltaf\eqnum*
$$
It is a homomorphism. The image is
 functions on $B$ and  operators acting on
such functions.

Given an arbitrary set of functions
$\{\phi_b(\ ), b\in B\}$, we can also apply difference operators
to the sufficies:
$$
\eqalign{
& {}_\de (ga)(\sum_{b\in B}c_b\phi_b(\ ))\ =
\ \sum_{b\in B}c_b g( q^b t^{-\rho})\phi_{b-a}(\ ),\
c_b\in \C.
}
\eqnu
\label\deltsuf\eqnum*
$$
It is an anti-homomorphism,
 i.e.
$$
{}_\de (GH)\ =\  {}_\de H\ {}_\de G \for \hbox{\ difference\
operators\ } G,H.
$$

{}From now on we will mostly use the {\it
renormalized Macdonald polynomials}
 $\pi_b(x)\equal p_b(x)/p_b(t^{-\rho})$ and their
discretizations:
$$ \pi_b(c)\equal {}^\de p_b (c)/{}^\de p_b (0) =
p_b (q^{c}t^{-\rho})/ p_b (t^{-\rho})=
\pi_c(b) \for b,c\in B_-.
$$

Given a symmetric polynomial $f\in \C[x]^W$, we
construct the operator $L_f$, go to its discretization
${}^\de L_f$,
and finally introduce the {\it recurrence operator}
 $\La_f= {}_\de L_f$  acting on the sufficies
$b\in B$ of any $\C$-valued
functions
$\phi_b (\ )$. We write $L_a,\La_a$ when $f$ is
 the monomial symmetric function $m_a, a\in B_-.$

\proclaim{Recurrence Theorem}
For arbitrary $a,b\in B_-, f\in \C[x]^W$,
$$
\eqalign{
&\La_{f}(\pi_b(x))\ =\  \bar{f}(x)\pi_b(x),\
\ \La_{a}(\pi_b(x)) =  \bar{m}_a(x) \pi_b(x),
}
\eqnu
\label\pbpb\eqnum*
$$
where  $\bar{f}(x)=f(x^{-1})$. The operators $\La$ (acting
on $b$)  do
 not produce  $\pi_c$ for $c\not\in B_-$.

\label\RECUR\theoremnum*
\endproclaim

{\it Proof.}
We can rewrite (\ref\Lf) as follows:
$$
\eqalignno{
&{}^\de L_f({}^\de p_b)\ =\ \bar{f}(q^{b}t^{-\rho})\ {}^\de p_b.
&\eqnu
\label\deLf\eqnum*}
$$
Replacing $p$ by $\pi$ and using the duality we yield:
$$
\eqalignno{
{}^\de & L_f(\pi_b(c))\ =\ \bar{f}(q^{b}t^{-\rho})
\ {}^\de \pi_b(c),\cr
& \La_f(\pi_b(c))\ =\ {}^\de\bar{f}(c)
\ {}^\de \pi_b(c),&\eqnu
\label\deLaf\eqnum*
}
$$
and (\ref\pbpb) if we can ensure that  $\La_f$ does not
create polynomials $p_c$ with the indices apart from $B_-$. The
latter can
be checked due to
 Lemma 2.4 from [C2]. When $a=-b_r, r\in O^*$
it results  directly from
Proposition \ref\LEAD.
We will give below another
proof independent of [C2].

\proofbox

To compare the theorem with the classical results about the products
of the characters (tensor products of the  irreducible
representations) take $t=q$. Then  corresponding Macdonald polynomials
can be obtained from the monomial functions by means of the shift operator
for $t=1$ which is very simple and  leads to the Weyl  formula
(this operator was introduced in [BG]).
Hence they are the characters of
the irreducible representations.
For minuscule  $-a$ in the case of $A_n$ (all dominant weights
are their algebraic combinations)
 the theorem was obtained by
Macdonald, Koornwinder, and then in [EK2].
 As to other root
systems (except $BC_1$ and certain special $t,q$), it seems to be new.

The theorem has many applications. For instance we can get another
prove of formula (\ref\EV) and moreover calculate the norms of $\pi_b$.
To demonstrate this let
$$
\eqalign{
& L_a\ =\ \sum_b b g_a^b(X),\ \La_a\ =\
\sum_b  {}_\de (g_a^b)\ {}_\de b,  \ \   a\in B_- ,\cr
 & L_a\ =\ \sum_b  f_a^b(X) b,\ \La_a\ =\
\sum_b {}_\de b\  {}_\de (f_a^b) ,\ w_0(a)\succeq b\succeq a.
}
\eqnu
\label\ella\eqnum*
$$
Then $g_a^b(q^c t^{-\rho})=f_a^b(q^{c+b} t^{-\rho})$.
The $w_0$-invariance of $L$ leads to the following relations:
$$
g_a^{w_0(b)}(t^{\rho})\ =\  g_a^b(t^{-\rho}),\ \
f_a^{w_0(b)}(t^{\rho})\ =\  f_a^b(t^{-\rho}).
\eqnu
\label\bprime\eqnum*
$$
We mention that  $\langle\mu_1\bar{m}_a\rangle= g_a^0(t^{\pm\rho})=
f_a^0(t^{\pm\rho})= \langle\mu_1{m}_a\rangle$
are the coefficients of the symmetrization
of $\mu_1$. Their  good description is one of the main
open problems in the Macdonald theory (we will consider it in
the next paper).

\proclaim{Proposition}
Setting $b^o\equal -w_0(b)$,
$$
\eqalignno{
&\sum_b f_a^b(t^\rho)\pi_{b^o}\ =\ m_{a^o}\ =\
\sum_b g_a^b(t^{-\rho})\pi_{b^o}\langle \pi_b,\pi_b\rangle^{-1}&\eqnu
\cr
\label\fma\eqnum*
&\langle \pi_b,m_a\rangle\ =\
g_a^b(t^{-\rho})\ =\ f_{a}^b (t^\rho)\langle \pi_b,\pi_b\rangle,\
f_b^b(t^\rho)\  =\ p_b(t^{-\rho}),  &\eqnu
\cr
\label\fqgq\eqnum*
&\langle\pi_b,p_b\rangle\ =\ \langle\pi_b,m_b\rangle\
=\ p_{b}(t^{-\rho})\langle \pi_b,\pi_b\rangle\ =
\ g_b^b(t^{-\rho})\ = \cr
&x_b(t^{-\rho})\prod_{a\in R_+^\vee}
\ \prod_{ 1\le j\le -(a^\vee,b)}
\Bigl(
{
1-q_a^{j}t_a^{-1}x_a(t^\rho)
 \over
1-q_a^{j}x_a(t^\rho)
}
\Bigr).
&\eqnu
\label\norms\eqnum*
}
$$
\label\FMG\theoremnum*
\endproclaim

{\it Proof.} One has:
$\langle m_{a^o},\pi_{b^o}\rangle=
 \langle \pi_{b}, m_a \rangle=
 \langle \mu_1\pi_{b} \bar{m}_{a} \rangle = $
 $\langle \mu_1\La_a\pi_{b} \rangle = $
$g_a^b(t^{-\rho})$. On the other hand,
$ {m}_{a^o}= \pi_0 {m}_{a^o}=
\sum_b f_a^b(t^\rho)\pi_{b^o}$. It gives all the relations
for $\{g,f\}$. The functions $g_b^b,f_b^b$ were calculated
in Proposition \ref\LEAD.
\proofbox

We note that (\ref\norms) can be directly deduced from
the formulas for  the norms of the Macdonald polynomials
$\{p_b\}$ from
[C2] adapted to  the present
paper. Indeed, given $b\in B_-$,
$$
\eqalign{
&\langle \mu p_b \bar{p}_b\rangle\ =\ \prod_{a\in R_+^\vee}(
(1-t_a x_a(t^\rho))/
(1-x_a(t^\rho)))\cr
&\prod_{a\in R_+^\vee, 0\le j< \infty}
\Bigl(
{
(1-q_a^{j+1}x_a(t^\rho q^{-b}))
(1- q_a^{j}x_a(t^\rho q^{-b}))
 \over
(1-t_a q_a^{j}x_a(t^\rho q^{-b}))
(1-t_a^{-1} q_a^{j+1}x_a(t^\rho q^{-b}))
}
\Bigr).
}
\eqnu
\label\NORM\eqnum*
$$
The necessary analitical continuation
is precisely as in [EK2] (the remark
after Theorem 3.9).
We must divide them by $\langle \mu\rangle$ to
go to the scalar product (\ref\innerpro):
$$
\eqalign{
&\langle  p_b, {p}_b\rangle\ =\ \prod_{a\in R_+^\vee}
\ \prod_{ 0\le j< -(a^\vee,b)}\cr
&\Bigl(
{
(1- q_a^{j+1}t_a^{-1}x_a(t^\rho))
(1- q_a^{j}t_a x_a(t^\rho))
 \over
(1- q_a^{j}x_a(t^\rho))
(1- q_a^{j+1}x_a(t^\rho))
}
\Bigr).
}
\eqnu
\label\scalar\eqnum*
$$
Then it is necessary just to divide
by
$$
\eqalignno{
&p_b(t^{-\rho})\ =\ x_b(t^{\rho})\prod_{a\in R_+^\vee}
\ \prod_{ 0\le j< -(a^\vee,b)}
\Bigl(
{
1- q_a^{j}t_a x_a(t^\rho)
 \over
1- q_a^{j}x_a(t^\rho)
}
\Bigr).
&\eqnu
\label\fbbev\eqnum*
}
$$

\vskip 10pt
 {\it Schwartz functions.}
We will use Macdonald's polynomials
to construct the eigenfunctions of the Fourier transform
which are also pairwise orthogonal
 with respect to the pairing
$[\![\ ,\ ]\!]$. The Gaussian comes naturally in this stuff
because of
the following theorem.

\proclaim{ Theorem}
 {i)} Adding   $ q^{1/2m}$, the following maps
can be uniquely extended to automorphisms of \HH\ , preserving each of
$T_1,\ldots,T_n,t$ and $ q$:
$$
\eqalign{
& \tau_+: \ X_b \to X_b,\ \ Y_r \to X_rY_r q^{-(b_r,b_r)/2},\
Y_\th \to X_0^{-1}T_0^{-2}Y_\th,      \cr
& \tau_-: \ Y_b \to Y_b,\ \ X_r \to Y_r X_r q^{(b_r,b_r)/2},\
X_\th \to T_0X_0Y_\th^{-1} T_0, \cr
& \om: Y_b \to X_b^{-1},\ X_r \to X_r^{-1}Y_r X_r q^{(b_r,b_r)},\
\ X_\th \to T_0^{-1}Y_\th^{-1}T_0, \cr
&\where b\in B,\ r\in O^*,\ X_0\ =   q X_\th^{-1}.
}
\eqnu
\label\tauom\eqnum*
$$
ii) The above maps give   automorphisms of \HH\  and
the  (elliptic
braid) group \BB\  generated by
the elements $\{X_b,Y_b,T_i,\pi_r, q^{1/2m}\}$ satisfying
the relations (i)-(vi) from Definition \ref\double  and
(\ref\Yb). Let  $\AA_o$ be the group  of its automorphisms modulo
the conjugations by the elements
from  the center $Z(\B)$ of the group
$\B$ generated by $\{T_1,\ldots,T_n\}$.
Considering the images of $\vep$ (see (\ref\vep)),$\tau_{\pm},\om$
in $\AA_o$ we obtain the homomorphism $GL_2(\Z)\to \AA_o$:
$$
\eqalign{
& \Bigl(\matrix  0 &-1\\ -1& 0\endmatrix \Bigr) \to \vep,\
  \Bigl(\matrix 1& 1\\ 0& 1 \endmatrix \Bigr) \to \tau_+,\cr
& \Bigl(\matrix 0& -1\\ 1& 0\endmatrix \Bigr) \to \om,\
\Bigl(\matrix 1& 0\\ 1& 1   \endmatrix    \Bigr) \to \tau_-.
}
\eqnu
\label\glz\eqnum*
$$
iii) The automorphism $\tau_-$ leaves the left ideal
$I$ from \ref\hat invariant and therefore acts in $\C_{ q,t}[x]$
identified with $\HH/I$.
\label\GLZ\theoremnum*
\endproclaim
{\it Proof.}
 The theorem can be deduced
from the topological interpretation of \BB\ from [C1] as
the fundamental group of the elliptic configuration space
(a proper  version of the product of $n$ copies of an elliptic
curve without the complexifications of the root hyperplanes and
divided by $W$). The standard action of the $GL_2(\Z)$ on the
periods of the elliptic curve results in the formulas  for
$\vep, \tau_+$. The remaining elements can be expressed in terms
of these two:
$$
\eqalign{
&\tau_-\  =\  \vep\tau_+\vep =  \vph\tau_+\vph,\ \ \om\  =\
\tau_+^{-1}\tau_-\tau_+^{-1} =
\tau_-\tau_+^{-1}\tau_-.
}
\eqnu
\label\taumin\eqnum*
$$
Since $GL_2(\Z)$ preserves the origin of the elliptic curve
the relations from this group are fulfilled up to conjugations
by the elements from $\B$. However the latter elements
must belong to
the center because $\tau_{\pm},\om$ fix the elements
$T_1,\ldots,T_n$.
 In the case of $A_n$,
the calculations are in fact  due to
J. Birman.

One can eliminate the topology and check
the statements  directly. Because we have already
(\ref\taumin) it suffices to calculate that
 $\om^4$ is the conjugation by
the element $T_{w_0}^2$, which follows from (\ref\tauom).
Topologically, a nontrivial step is to see that
$\tau_+$ (which is evidently a  pullback of
$\Bigl(\matrix 1&1\\ 0&1\endmatrix\Bigr) $)
is really an automorphism of \BB\ (for $A_n$ it is not difficult).
Fortunately we can
make it algebraic as well.

Setting
$x_b= q^{z_b}, \ z_{a+b}=z_a+z_b, \ z_i=z_{b_i},\
a(z_b)= z_b-(a,b), \ a,b\in \R^n,$
we introduce the {\it Gaussian function}
$\ga\ =\  q^{\Sigma_{i=1}^n z_i z_{\al_i}/2}$,
which can be considered as a  formal series in $x_i-1, \log q$
and satisfies the following (defining) difference relations:
$$
\eqalign{
&b_j(\ga)\ =\  q^{(1/2)\Sigma_{i=1}^n (z_i-(b_j,b_i))
(z_{\al_i}- q_i^j)}\ =\cr
& \ga  q^{-z_j+ (b_j,b_j)/2 }\ =\  x_j^{-1}\ga  q^{(b_j,b_j)/2}
\for 1\le j\le n.
}
\eqnu
\label\gauss\eqnum*
$$

The Gaussian function commutes with $T_j \for 1\le j\le n$
because it is $W$-invariant.
When $b_r$ are minuscule ($r\in O^*$), we  use directly  formulas
(\ref\Yb, \ref\Demaz) to check that
$$
\ga(X)Y_r\ga(X)^{-1}\ =\  X_r q^{-(b_r,b_r)/2} Y_r \ =
\ \tau_+(Y_r).
$$
A straightforward calculation gives that
$$
\eqalign{
&\ga(X)T_0\ga(X)^{-1}\ =\  \tau_+(T_0) = X_0^{-1}T_0^{-1},\cr
&\tau_-(T_0)\ =\  T_0,\
\om(T_0)\ =\ X_\th^{-1}Y_\th^{-1}T_0,
}
\eqnu
\label\gato\eqnum*
$$
which yields (\ref\tauom).

Here the multiplication by $\ga$ is not well defined
(only the conjugation). It preserves the bigger space
of meromorphic functions in $z_b$ which is
an \HH\ - module as well.  In $\C_{ q,t}[x]$ the
picture is  opposite. The automorphism $\tau_-$ is inner since
it fixes $\{Y, T\}$.
\proofbox
We note that the automorphisms $\{\vep,\tau_{\pm},\om\}$
are {\it projectively unitary} with respect to the Macdonald
pairing. Indeed, they send the unitary operators
$\{X,Y,T, q,t\}$ generating \HH\ to  unitary ones.

The next claim is that the
{\it Schwartz functions} $\{\ga^{-1} p_b,\ b\in B_-\}$ defined
for the Mac\-do\-nald
poly\-no\-mials $\{p_b\}$ are eigenfunctions of the Fourier transform
and  pairwise or\-tho\-go\-nal with respect to
the Fo\-urier pairing $[\![\ ,\ ]\!]$. We will reformulate and
prove it algebraically
representing Schwartz functions
as eigenfunctions of self-adjoint operators
with respect to  $\vph,\vep$.

\proclaim {Proposition}
The $W$-invariant operators $\l_f^\ga\equal \ga^{-1} \l_f\ga$ defined for
 $f\in \C[x]^W$ obey the following properties
(see Theorem \ref\LF):
$$
\eqalign{
&\l_f^\ga (\ga^{-1} p_b)\ =\ f(t^\rho q^{-b})(\ga^{-1} p_b),\cr
&\vph(\l_f^\ga)\ =\ \l_f^\ga, \ \vep(\l_f^\ga)\ =\ \l_{\bar{f}}^\ga.
}
$$
The
corresponding eigenvalues (for all $f$)
distinguish different $\ga^{-1} p_b$.
\endproclaim
\label\FEIGEN\theoremnum*
\proofbox

There is another way to make the statement about
the action of the Fourier transform on the Schwartz functions
algebraic. The following formula gives this and moreover the
exact eigenvalues (up to proportionality):
$$
\eqalign{
& L_{p_b}(\ga^{-1})\ =\  q^{-(b,b)/2}
x_b(t^\rho) {p_b}\ga^{-1}, \ b\in B_-.
}
\eqnu
\label\Lfga\eqnum*
$$
The coefficents of proportionality are calculted by means
of Proposition \ref\LEAD (see also the end of the paper).

One
 can introduce formally the {\it Fourier transform} $F( f)$ of
a function $f$ by means of the relation (see (\ref\innerpro)
 and (\ref\Fourier)):
$$
\eqalign{
& \langle g, F(f)\rangle\ =\ [\![g,f]\!]\  =\
\l_{\bar{g}}(f)(t^{-\rho})
}
\eqnu
\label\ftran\eqnum*
$$
 valid for all polynomials $g$ (or even for more general classes
of functions).
The associated  transformation  on the
operators is
exactly $\vep$ from (\ref\vep).
When $t=1$ it is the classical Fourier transform:
$$F( f)(\la)\ =\ \hbox{const}\int_{\R^n}  q^{\Si z_i\la_{\al_i}}
{f}^+(z) dz_1\ldots dz_n,
 $$
where ${}^+$ is the formal complex conjugation $ q^+= q^{-1},\
t^+=t^{-1}$. It is the standard conjugation if $ q,t$ belong
to the unit circle.
The appearance of ${f}^+$ instead of $f$
 makes the transform involutive (and is almost inevitable
in the difference theory).

Hopefully  $\tau_{\pm}$ become inner
and the above claims  about Schwartz functions
 can be made rigorous after the
discretization (see the beginning of the section).
When  $ q$
is a root of unity the  procedure is completely
algebraic.

%
%
%
%
\section {Roots of unity}

Let us assume that $ q$ is a primitive
$N$-th root of unity for $N\in \N$
and  consider $t$ as an indeterminate parameter.
More precisely, we will operate over the field
$\Q_t^0\equal \Q( q_0,t)$ where we fix $ q_0$
such that $ q_0^{2m}= q$ ($ q_0$ belongs to a proper
extension of $\Q$). Actually all formulas will  hold
even over the localization of $\Z[ q_0,t]$
by $t^{r_1} q^{s_1}(1-t^{r_2} q^{s_2})\neq 0$, $r_i,s_i\in \Z$.
The pairing
$$
\eqalign{
B\times B\ni a\times b\to  q^{(a,b)}\equal  q_0^{2m(a,b)}
}
\eqnu
\label\BtimeB\eqnum*
$$
acts through $B_N\times B_N$, where
 $B_N\equal B/K_N$,  $K_N$ is its radical.

Following  the previous section
we restrict the functions $\{m_b\}$ and
the operators $\{L_b=L_{m_b}\}$ to the
lattice $B$ using the pairing (\ref\BtimeB).
We need to consider the quotient $B_N$ only.
The $L$-operators
are well defined over $\Q_t^0$ since their
 denominators
are  products of the binomials $(x_a q^k-1)$
for $a\in R^\vee, k\in \Z$. The latter remain
non-zero when evaluated at $ q^b t^{-\rho}$
since $(a,\rho)\neq 0$ ($x_a(t^\rho)$
always contain $t$).
More exact information
about the properties of these coefficients can
be extracted from Proposition \ref\FMG.

Let $B_-(N)$ be a fundamental domain of the
group $K_N$. It means that the map
$B_-(N)\to B_N$
is an isomorphism.
Further we identify these two sets, putting
$$
B_-(N)\ =\
\{\be^1,\ldots, \be^d \}\  =\  B_N, \where d=|B_N|,\ \be^1=0,
\eqnu
\label\bei\eqnum*
$$
and denote the image of $b\in B$
in $B_N$ by $b'$.
We note that
 ${}^\de m_a (c) = m_a(t^{-\rho} q^{c})$ for $ a\in A$
separate $\{\be^i\}$. It is true since
 $w(t^\rho q^{-b})\neq t^\rho q^{-c}$
for any $w\in W$ if  $b'\neq c'$
(we remind that $t$ is generic). We assume that
$-w_0(B_-(N))=B_-(N)$ for the longest element $w_0$.

Let us consider (temporarily) the case when
$N$ is coprime with the order $ |B/A|= |O|$ taking
$ q_0^2= q^{1/m}$ in the $N$-th roots of unity.
Then $K_N=  NP\cap B$ for the weight lattice
$P=\oplus_{i=1}^n\Z\om_i$ generated by the  $\om_i$
(dual to $a_i$). We can take the following fundamental domain
$$
\eqalign{
 &B_-(N)\ =\ \{b= -\sum_{i=1}^n k_ib_i\in B_-\} \hbox{\ such\ that}\cr
&0\le k_i< N \if (2/\nu_i, N)=1,\ 0\le k_i<\nu_i N/2
\hbox{\ \ otherwise}.
}
\eqnu
\label\kineq\eqnum*
$$
We remind that $2/\nu_i=2/(\al_i,\al_i)= 1,2,3$
(see Introduction).

Let us demonstrate that the Macdonald polynomials $p_b$
are well defined for $b\in B_-(N)$ (later we will
see that they  always exist). We  introduce them
 directly from (\ref\Lf), using that the $L$-operators
preserve any subspaces
$$
U_b\ =\ \oplus_{ c\succ b} \Q_t^0 m_c,
\and  U_b\oplus \Q_t^0 m_b,\ b,c\in B_-.
$$
Here it is necessary to check that given $B_-\ni c\succ b$,
there exists at least one $ a\in B_-$ such that
$m_a(t^\rho q^{-b})\neq  m_a(t^\rho q^{-c})$.
Then the eigenvalues of $L$ will separate $p_b$ from
the elements from $U_b$ and we can argue by induction.
If $A_+\ni c-b\in NP$ then there is $a_i (1\le i\le n)$
such that $(c-b,a_i)>0$ (the form $(\ ,\ )$ is positive).
Since $(A,A)\subset (2/\nu_i)\Z$,
$$
 (-b,a_i)\ge -(b,a_i)+(c,a_i)
\ge (2N/\nu_i)/(2/\nu_i,N),
$$
which contradicts (\ref\kineq).

We note that the norms $\langle p_{b},p_{b}\rangle$
are non-zero for all these polynomials
which results from (\ref\scalar). The same holds
for $\langle \pi_{b'},\pi_{b'}\rangle,\ b'\in B_-(N),$
due to (\ref\EV).

 The discretizations of
$\pi_{b'}$ are well defined too and
$\pi_{b'}(c)$ depends only on the image $c'$ because
$\pi_{b'}$ is a linear combinations of $m_a, a\in B_-$.
Finally,  the
 $\{\pi_{\be^i}(c')\}$ form a basis in the
space $V_N\equal \hbox{Funct}(B_N, \Q_t^0)$ of all $ \Q_t^0$-valued
 functions on $B_N$. Indeed they are non-zero and
the action of the $\{{}^\de L_a\}$ ensures that they are linearly
independent.

{\it The end of the  proof of Theorem \ref\RECUR }.
First of all, let us rewrite formally relation
 (\ref\pbpb) for $f=m_a$
 as follows:
$$
\eqalign{
&\bar{m}_a(x)\pi_b(x)\ =\ \La_a^{-} (\pi_b(x))+
 \sum_{e\not\in B_-} M_{ab}^{e} \pi_{e}(x).
}
\eqnu
\label\mpie\eqnum*
$$
Here $e$ form a finite set $E=E(a,b)\ (E\cap B_-=\emptyset)$,
$M_{ab}^e$ are rational functions of $ q,t$.
The
truncation $\La_a^-$
of $\La_a$ is uniquelly determined by the
 condition that it does not contain the shifts moving $b$ to
elements apart from $B_-$. Assuming that $N$ is sufficiently
big ($B_-(N)$ must contain $b, b+a$ and $c\in B_-$ such that
$ c\succeq b-w_0(a)$)
the discretization gives the relation (see (\ref\deLaf)):
$$
\eqalign{
&\bar{m}_a(c)\pi_b(c)\ =\ \La_a^{-} (\pi_b(c))+
 \sum_{e\not\in B_-} M_{ab}^{e} \pi_{e}(c),\  c\in B_-(N),
}
\eqnu
\label\mpiede\eqnum*
$$
for $ m_a(c)={}^\de m_a (c)$.
Here   $\pi_e=\pi_{e'}$ for $e\in E$. This substitution was
impossible before the discretization. We remind
that  the  formula  with
 $\pi_c(e)$ in place of  $\pi_e(c)$ is always true.
Because $c$ is taken
from $B_-(N)$ the  discretization of $\pi_c$
 exists. Therefore  we can go from $e$ to $e'$,
 and then replace  $\pi_c(e')$ by $\pi_{e'}(c)$.
As to $M_{ab}^e$, they  are  the values of the coefficients
of ${}^\de L_a$ and are also well-defined when
$ q^N=1$  (enlarging $N$ we can get rid of singularities
in $ q$ even if $M$ are arbitrary rational).

On the other hand :
$$
\eqalign{
&\bar{m}_a(x)\pi_b(x)\ =\
 \sum_{e\in B_-} K_{ab}^{e} \pi_{e}(x),
}
\eqnu
\label\mpix\eqnum*
$$
where the coefficients $K_{ab}^e$ are rational functions
of $ q,t$,
$\{e\}$ form a finite set $E_- =E_-(a,b)$. The discretization
 gives that
$$
\eqalign{
&\bar{m}_a(c)\pi_b(c)\ =\
 \sum_{e\in B_-} K_{ab}^{e} \pi_{e}(c), \ c\in B_-(N).
}
\eqnu
\label\mpixde\eqnum*
$$
We pick  $N$ to avoid possible singularities.

Since $N$ is sufficiently big, the eigenvalues of the
$L$-operators distingwish all $\pi_{e'}(c)$  for $e\in
E\cup E_-$. It holds only for generic $t$ (say, when
$t=1$ it is wrong). Comparing (\ref\mpiede) and
(\ref\mpixde) we conclude that $M_{ab}^e=0$ for all $e\in E$,
when $ q^N=1$.
Using again that $N$ is arbitrary (big enough, coprime with $|O|$)
we get
that the actions
of $\La_a^-$ and $\La_a$   coincide on $\pi_b$, i.e.
the latter operator does not create the indices  not
from $B_-$.
\proofbox

Let us go back to the general case (we drop the condition
$(N,|O|)=1$).
Once the Recurrence Theorem has been established we can
use Proposition \ref\FMG without any reservation.
It readily gives that the Macdonald polynomials $p_b$
and  ${}^\de\pi_b$ are well defined for
arbitrary $b\in B_-$ because
 $f_b^b(t^{-\rho})\neq 0$ (Proposition
\ref\LEAD). Moreover,
${}^\de\pi_b = {}^\de\pi_c$ if and  only  if
$b'=c'$,
 and the {\it restricted Macdonald polynomials} $\pi_{\be^i}(c'),
1\le i\le d$ (see (\ref\bei)) form a basis in
 $V_N= \hbox{Funct}(B_N, \Q_t^0)$.
Indeed,   $\pi_{\be^i}(c')$ are eigenvectors of the
 ${}^\de L$-operators separated by the eigenvalues. They
 are always non-zero since
$\pi_b(0)=1$. Hence they are linearly independent over
$\Q_t^0$ and form a basis in $V_N$. Every ${}^\de\pi_b$
is an $L$-eigenvector and coincides
with  one of them (when $\be^i=b'$).
 Similarly,
${}^\de L_{\pi_b} = {}^\de L_{\pi_c}$ if and  only  if
$b'=c'$ because the latter condition is necessary and
sufficient to ensure the coincidence of the sets of
eigenvalues.
We will also use the basis of
the {\it delta-functions}
$\de_{\be^i} (\be^j)\equal\de_{ij}$ separated
by the action of $\{{}^\de m_a\}$.

Let $\a$ be the algebra of the elements of \HH\
 commuting with $\{T_1,\ldots,T_n\}$.
All the (anti-)automorphisms under consideration preserve
it.

\proclaim {Proposition}
 The discretization map
$ \a \ni A \to A_\dagger= [\hat{A}]_{\dagger}\to
{}^\de A_\dagger
$
 supplies
$V_N$ with the structure of
an ${}^\de \a$-module which is irreducible.
The algebra ${}^\de \a $ is generated by the discretizations of
$L_a, m_a$ for   $a\in A$. The Fourier pairing is well defined
on $V_N$ and induces the anti-involution
$\vph$ ($\bar{m}_a\to L_a\to \bar{m}_a$).
\endproclaim
\label\IRRED\theoremnum*
{\it Proof.}
 The radical of the Fourier pairing (\ref\Fourier)
 contains the kernel of the discretization
map $\Q_t^0[x]^W\to V_N$. Its restriction
 to $V_N$ is  non-degenerate since
 $$
\eqalign{
\Pi=(\pi_{ij}), \where \pi_{ij}\ =\
[\![\pi_{\be^i},
 \pi_{\be^j}]\!]= \pi_{\be^i}(\be^j),
}
\eqnu
\label\pimat\eqnum*
$$
is the matrix connecting
the bases $\{\pi\}$ and $\{\de\}$.
The coresponding anti-involution $\vph$
transposes $\bar{m}_a$ and
$L_a$.
Thus  $V_N$
 is semi-simple.

If $V_N$ is reducible then
$\pi_{\be^1}=1$   generates a
proper $\a$-submodule ($\neq V_N$). But it  takes
 non-zero values at any points of $B_N$. Hence
its $\a$-span must contain all $\de_{\be^i}$. We come to
a contradiction.
\proofbox

{\it When $t$ are special.} Till the end
of the paper  $t_\nu= q_\nu^{k_\nu}$ for
$k_\nu\in \Z_+,\  \nu\in \nu_R$. The $L$-operators
act in  $\Q^0[x]$ for $\Q^0\equal \Q( q_0)$.
Let $J\subset \Q^0[x]^W$ be the radical of the pairing
$[\![\ ,\ ]\!]$. It is an ideal and an
$\a$-submodule. The quotient (a ring and an $\a$-module)
$\v=\Q^0[x]^W/J$ is finite dimensional over $\Q^0$.
It results from Proposition \ref\IRRED as generic
$t$ approaches $ q^k$.

The set $\tB$ of maximal ideals
of $\v$ can be considered as a subset of
any fundamental domain
(in $B_-$) with
respect to the action of the
groups $K_N$ and $\ W\ni w: b\ \to\ $ $w(b+\sum_\nu k_\nu\rho_\nu)
-$ $\sum_\nu k_\nu\rho_\nu$.
We put
$$
\tB\ =\
\{\tbe^1,\ldots, \tbe^\pa \}\  \subset\  B_-, \ \tbe^1=0,
\eqnu
\label\tbei\eqnum*
$$
assuming that $-w_0(\tB)=\tB$ for the longest element $w_0$.
We remind that the point $0$
corresponds to the evaluation at $t^{-\rho}$.

By the construction,   all $\{L_a\}$-eigenvectors
in $\v$ have pairwise distinct eigenvalues (the difference
of any two of them  with the same sets of
eigenvalues and coinciding evaluations at $t^{-\rho}$
belongs to $J$). Applying this to  $\pi_0=1$ generating
$\v$  as an $\a$-module
we establish the irreducibility of $\v$ (use the pairing
$[\![\ ,\ ]\!]$ and follow Proposition \ref\IRRED).

Next, we will introduce the {\it restricted Macdonald
pairing} (cf. (\ref\mu),(\ref\innerpro)):
$$
\eqalign{
\langle f(x),g(x)\rangle'&\equal\sum_{c\in B_N}
\mu'(c)f(c)\bar{g}(c) \for f,g\in \Q^0[x]^W,\cr
\mu'&\ =\ \prod_{a\in R^\vee_+}
\prod_{i=0}^{k_a-1}(1-x_a q_a^i)(1-x_a^{-1} q_a^{i}).
}
\eqnu
\label\resinner\eqnum*
$$
Here $B_N=B/K_N,\ \mu'(c)={}^\de\mu'(c)=\mu'(t^{-\rho} q^c)$.
The usage of symmetric $\mu'$ in place of $\mu$ (which does not
alter $\langle\ ,\ \rangle$ up to proportionality
 for generic $t$) is somewhat more
convenient for roots of unity.

The same verification as in [C2],  Proposition 4.2
gives that
$$
\langle Af,g\rangle ' \ =\ \langle f,A^*g\rangle',\
A\in [\hat{\a}]_\dagger,
\eqnu
\label\resinv\eqnum*
$$
for the anti-involution ${}^*$ from (\ref\star)
considered on  the operators
 from $\a$ acting on symmetric polynomials.
The restricted norms can be calculated by means of
the same shift operators till the latter
do not vanish (see [C2], Corrolary 5.3), which leads to
the   formulas (\ref\NORM) up to a common coefficient of
proportionality.

Let us assume that:
$$
 q_a^{(\rho_k,a)+i}\neq 1 \hbox{\ for\ all\ } a\in R^\vee_+,\
i=-k_a+1,\ldots,k_a-1,
\eqnu
\label\krho\eqnum*
$$
where $\rho_k= \sum_\nu k_\nu\rho_\nu$.
This restriction makes
the construction below non-empty. In the simply
laced case ($A,D,E$), it is equivalent to the condition
 $N\ge
k((\rho,\th) + 1)$.

\proclaim {Lemma}
The natural map $\v\to \tV \equal \hbox{Funct}(\tB,\Q^0)$
is an isomorphism
which supplies $\tV$ with the structure of a non-zero
irreducible $\a$-module. Both
 pairings $[\![\ ,\ ]\!], \langle\ ,\ \rangle$ are
well defined and non-degenerate on $\tV$.
\endproclaim
\label\COIN\theoremnum*

{\it Proof.}
The radical $J'$ of the pairing
$\langle\ ,\ \rangle'$ in $\Q^0[x]^W$ is an $\a$-submodule.
It equals the space of all functions $f(x)$ such that
${}^\de f(c)=0$ for $c$ from the subset $B'\subset \tB_-(N)$
where ${}^\de\mu'$ is non-zero. The  sets  $\tB$ and $B'$  contain
$0$. Indeed, $[\![1,f]\!]=f(0)$ and
$\mu'(t^{-\rho})\neq 0$ because of  condition (\ref\krho).
Hence the linear span $J+J'$ (that is an $\a$-submodule)
does not coincide with the entire $\Q^0[x]^W$, and
the  irreducibility of $\v$ results in $J+J'=J$.
\proofbox

Introducing now the delta-functions $\tde_i=\de_{\tbe^i}$,
we can define the $\pi$-functions $\{\tpi_i\}$
from the orthogonality and evaluation conditions
$$
[\![\tpi_i,\tde^j]\!]= C_i\de_{ij} \and \tpi_i(0)=1,\
, 1\le i,j\le \pa.
\eqnu
\label\tpitde\eqnum*
$$
They are eigenvectors of the $L_a$-operators with
the eigenvalues $\bar{m}_a(\tbe^i)$ and linearly
 generate $\tV$. The sets of eigenvalues are pairwise
distinct and $\langle \tpi_i,\tpi_i\rangle'\neq 0$.

Actually the $\pi$-functions are the discretizations of
certain restricted Macdonald's polynomials and the
above scalar products can be calculated explicitly
but we will not discuss this here.
Informally,  $\tB$ is a collection of
 all weights (up to $K_N$)
corresponding to the Macdonald
polynomials of non-zero $q,t$-dimension.
Amyway, the restrictions to $\tB$ of any well defined
polynomial $\pi_b$ coincide with one of $\tpi_i$.

We will also use that
$(a,r_\nu)\in N\Z, $ where $r_\nu=(2/\nu)\rho_\nu\in B$,
and impose one more restriction:
 $$
q^{(a,a)/2}=q_0^{m(a,a)}\ =\ 1   \for a\in K_N, \nu\in\nu_R.
\eqnu
\label\aarho\eqnum*
$$
If $ q_0$ is a primitive root of degree $2mN$ then
$K_N=NQ\cap B$ for the root lattice $Q=\oplus_{i=1}^n \Z\al_i$
(see (\ref\BtimeB)). This condition  obviously holds true
for even $N$ (all roots systems). For odd $N$, it is necessary
to exclude  $B_n, C_{4l+2}$. In the latter case, $B\subset Q, m=1$
and we can
pick $q_0$ in the roots of unity of degree $N$.

Let us equip  the matrix algebra
$\hbox{End}_{\Q^0}(\tV)$  with the complex
conjugation ${}^+:   q\to  q^{-1}$
acting on the entries and   the hermitian
transposition  ${}^\dagger$ (the composition of the
matrix transposition ${tr}$ and ${}^+$).

\proclaim {Theorem}  Introducing $\tPi=(\tpi_{i}(\tbe^j))$
(see (\ref\pimat,\ref\tbei)) for $\tpi_i=\pi_{\tbe^i}$,
 let
$$
\eqalign{
&\t_+ = \hbox{Diag}( q^{(\tbe^i,\tbe^i)/2}x_{\tbe^i}(t^{-\rho})),\
\t_- = \Pi\t_+^{-1}\Pi^{-1},\ \Om  =
\t_+^{-1}\t_-\t_+^{-1}.
}
\eqnu
\label\DETOm\eqnum*
$$
Then the following map gives
 a non-zero projective representation of the group $SL_2(\Z)$:
$$
\Bigl(\matrix 1&1\\  0&1\endmatrix\Bigr) \to \t_+,\
\Bigl(\matrix 1&0\\  1&1\endmatrix\Bigr) \to \t_-,\
\Bigl(\matrix 0&-1\\  1&0\endmatrix\Bigr)\to \Om.
$$
Moreover,
 $\Pi^+ = \w_0 \Pi = \Pi\w_0,\
\t_+^{\dagger} = \t_+^{-1},\
\t_-^{\dagger} = \t_-^{-1},\
\Om^{\dagger} = \Om^{-1},
$
where $\w_0= (w_{ij}),\ w_{ij}= \de_{ij^o},\
\tbe^{j^o} = -w_0(\tbe^j)$ for the longest element $w_0$.
\label\ETD\theoremnum*
\endproclaim

{\it Proof.} We will start with the action of the
${}^+$ on  $\tpi_i(c)$  and the delta-functions:
$$
\eqalign{
&(\tbe^i)^+\ =\ \tbe^{i^o}= -w_0(\tbe^i),\ \tde_{i}^+\
=\ \tde_{{i^o}},\ \tpi_{i}^+\ =\ \tpi_{i}.
}
\eqnu
\label\plusV\eqnum*
$$
Given an arbitrary $f\in \tV$,  $f(\tbe^i)^+=
f^+(\tbe^{i^o})$. Together with the relation
$\tPi^{\hbox{tr}}=\tPi$ (the duality) it gives
the formulas with ${}^\dagger$.

 The multiplication by $\ga$ acts naturally
in $\tV$. The formula for $\t_+$ describes it
(up to a constant factor) in the basis of delta-functions.
 Really,
$$
\eqalign{
&\ga(c)\ =\ {}^\de \ga (c)\ =\
  q^{\Sigma_{i=1}^n \ze_i \ze_{\al_i}/2}\for \cr
&\ze_i\  =\ \log_q( x_i( q^c t^{-\rho}))\ =\
(b_i,c)+\log_{q} (\prod_\nu t_\nu^{-(\rho_\nu,b_i)}).
}
\eqnu
\label\taugau\eqnum*
$$
Hence $ \ga(c) = g  q^{(c,c)/2}x_c(t^{-\rho})$ for
$g =  q^{(\rho_k,\rho_k)}.$
Since  the  matrix $\t_+$  is important up to proportionality
one can  drop the constant $g$.
We see that changing $c$ by any elements from $K_N$
does not influence $\ga(c)$ because of the condition
(\ref\aarho), which makes the multiplication by $\ga$
well defined.

Next, the automorphism $\tau_- = \vph \tau_+\vph$
corresponds to $\tPi \t_+^{-1}\tPi^{-1}$, and the matrix $\Om$
from (\ref\DETOm) induces
 $\om=$ $\tau_+^{-1}\vph\tau_+\vph\tau_+^{-1}$
 in the same delta-basis.
Indeed, $\tau_-$ is the application
of $\ga(Y)$. It multiplies $\pi_b$ by $\ga(t^\rho q^{-b})$
whereas $\ga(X)$ multiplies $\de_b$ by $\ga(t^{-\rho} q^{b})$
(so we need to  inverse  $\t_+$). More formally,
one can use the equation
$[\![\t_+^{-1}f,g]\!]=[\![f,\t_-g]\!]$.

Finally, any relations from $SL_2(\Z)$ hold for these matrices up
to proportionality. It results directly from
 Theorem \ref\GLZ and the irreducibility
of $\tV$.
\proofbox

Actually we have come even to a stronger statement.

\proclaim {Corrolary}
The above map can be extended to a projective
representation
of $GL_2(\Z)$:
$$
\eqalign{
&\Bigl(\matrix 0&-1\\ -1&0\endmatrix\Bigr)\to \e
\equal \tPi\d^{-1}\tPi^+\si\tPi^{-1}\ =\
 \tPi\d^{-1},
}
\eqnu
\label\glproj\eqnum*
$$
where $ \d  =  \hbox{Diag}(\langle
\tpi_{i},\tpi_{i}\rangle'),\
\e$ belongs to the matrix algebra
$\hbox{End}_{\Q^0}(\tV)$ extended by the complex conjugation
$^+$ denoted by $\si$. Modulo proportionality,
$\Om$ equals $\d\tPi^{-1}$.
\endproclaim
\label\GLPROJ\theoremnum*

{\it Proof.}
We can introduce  the Fourier transform
 (\ref\ftran), giving $\vep$ on the operators:
$$
\eqalign{
&F({\tpi}_{j})\ =\ \sum_{i=1}^\pa \tpi_{i}(\tbe^j)^+
 \tpi_{i}/\langle \tpi_{i},\tpi_{i}\rangle'.
}
\eqnu
\label\foures\eqnum*
$$
Then $F$
is written  in the basis $\{\tpi\}$ as
$\d^{-1}\tPi^+\si$.
In the basis $\{\tde_{i}\}$, it will be exactly $\e$
from (\ref\glproj).
Here we applied the formula
$$(\tPi \d^{-1}\si)^2\ =\
\tPi^+\d^{-1}\tPi\d^{-1}\ =\ \hbox {const }
$$
resulting
 directly from the relation $\vep^2=1$ and the irreducibility
of $\tV$. Adding $\e$ to  (\ref\DETOm)
we get a projective representation
of $GL_2(\Z)$.

Let us check that
$\Om$ is proportional to $\d\Pi^{-1}$.
We use that the Fourier transform is diagonal in the
basis $\ga^{-1}\tpi_{i}$ (which
follows from Proposition \ref\FEIGEN).
It means that $\tau_+ F\tau_+^{-1}$ is diagonal
in the basis $\{\tpi\}$ and gives that
$$
\eqalign{
&(\tPi^{-1}\t_+\tPi)(\d^{-1}\tPi^+)(\tPi^{-1}\t_+^{-1}\tPi)^+\ =\
\f\equal  \hbox{Diag}(\phi_{\tbe^i})\for\cr
&   F(p_b\ga^{-1}) = \phi_b p_b\ga^{-1}\
 = \phi_0 L_{p_b}(\ga^{-1}),\
 \phi_b/\phi_0\ =\ x_b(t^\rho) q^{-(b,b)/2}.
}
\eqnu
\label\Fdiag\eqnum*
$$
The first expression is proportional to
$\tPi^{-1}\t_+\tPi\t_+\tPi^{-1}\d$ which is
$(\Om\tPi)^{-1}\t_+^{-1}\d$.
Then we observe that  $\t_+\f$ is a constant matrix.
Here the coefficients $\phi_b/\phi_0$ can be easily calculated
(once we have the proportionality). By the way,
taking rather big $N,k$ we come
to the proof of
formula (\ref\Lfga).
\endproof

The theorem
generalizes the construction due to Kirillov [Ki]
(in the case of $A_n$).
His approach is different and
 based on quantum groups at roots of unity.
As to arbitrary roots, the simplest case of our theorem
when $t= q$ is directly related to
Theorem 13.8 from [K].
The weights $\{\tbe^j\}$ exactly
correspond to  the representations of non-zero $ q$-dimension.
The Macdonald polynomials are the characters and
always exist  in this case.

We expect that the above considerations
can be applied with
proper changes to generic
$ q,t$ in the analitic setting
and hope that they could be useful to renew
elliptic functions towards the Ramanujan theories.

%
%
%
%
%
%
\AuthorRefNames [BGG]
\references
\vfil

[BG]
\name{F.A. Berezin}, and \name {I.M. Gelfand},
{Some remarks on the theory of spherical functions on
symmetric Riemannian manifolds}, Transl.,
II.Ser. Am. Math. Soc. {21} (1962), 193--238.

[B]
\name{N. Bourbaki},
{ Groupes et alg\`ebres de Lie}, Ch. {\bf 4--6},
Hermann, Paris (1969).

[C1]
\name{I. Cherednik},
{  Double affine Hecke algebras,
Knizhnik- Za\-mo\-lod\-chi\-kov equa\-tions, and Mac\-do\-nald's
ope\-ra\-tors},
IMRN (Duke M.J.) {  9} (1992), 171--180.

[C2]
\bibline,{ Double affine Hecke algebras and  Macdonald's
conjectures},
Annals of Mathematics {141} (1995), 191-216.

[C3]
\bibline,
{ Induced representations of  double affine Hecke algebras and
applications},
Math. Research Letters { 1} (1994), 319--337.

[C4]
\bibline,{ Difference-elliptic operators and root systems},
IMRN {1} (1995), 43--59.

[C5]
\bibline,
{ Integration of quantum many- body problems by affine
Knizhnik--Za\-mo\-lod\-chi\-kov equations},
Pre\-print RIMS--{  776} (1991),
(Advances in Math. (1994)).

[C6]
\bibline,
{ Elliptic quantum many-body problem and double affine
 Knizhnik - Zamolodchikov equation},
Commun. Math. Phys. (1995).

[D]
\name{C.F. Dunkl},
{ Hankel transforms associated to finite reflection groups},
Contemp. Math. {138} (1992), 123--138.

[EK1]
\name {P.I. Etingof}, and \name {A.A. Kirillov, Jr.},
{Macdonald's polynomials and representations of quantum
groups}, Math.Res.Let. {1} (1994), 279--296.

[EK2]
\bibline,
{Representation-theoretic proof of the inner product and
symmetry identities for Macdonald's polynomials},
Compositio Mathematica (1995).

[GH]
\name{A.M. Garsia}, and \name{M. Haiman},
{A graded representation model for Macdonald's polynomials},
Proc. Nat. Acad. Sci. USA {90}, 3607--3610.

[J]
\name{M.F.E. de Jeu},
{The Dunkl transform }, Invent. Math. {113} (1993), 147--162.

[H]
\name {S. Helgason},
{Groups and geometric analysis},
Academic Press, New York (1984).

[He]
\name{G.J. Heckman},
{  An elementary approach to the hypergeometric shift operators of
Opdam}, Invent.Math. {  103} (1991), 341--350.
Comp. Math. {  64} (1987), 329--352.
v.Math.{  70} (1988), 156--236.

[K]
\name {V.G. Kac},
{Infinite dimensional Lie algebras},
Cambridge University Press, Cambridge (1990).

[Ki]
\name {A. Kirillov, Jr.},
{Inner product on conformal blocks and Macdonald's
polynomials at roots of unity}, Preprint (1995).

[KL1]
\name{D. Kazhdan}, and \name{ G. Lusztig},
{  Proof of the Deligne-Langlands conjecture for Hecke algebras},
Invent.Math. {  87}(1987), 153--215.

[KL2]
\bibline,
{ Tensor structures arising from affine Lie algebras. III},
J. of AMS {  7}(1994), 335-381.

[KK]
\name{B. Kostant}, and \name{ S. Kumar},
{  T-Equivariant K-theory of generalized flag varieties,}
J. Diff. Geometry{  32}(1990), 549--603.

[M1]
\name{I.G. Macdonald}, {  A new class of symmetric functions },
Publ.I.R.M.A., Strasbourg, Actes 20-e Seminaire Lotharingen,
(1988), 131--171 .

[M2]
\bibline, {  Orthogonal polynomials associated with root
systems},Preprint(1988).

[M3]
\bibline, {  Some conjectures for root systems},
SIAM J.Math. Anal.{  13}:6 (1982), 988--1007.

[MS]
\name {G. Moore}, and \name{N. Seiberg},
{Classical and quantum conformal field theory},
Commun. Math. Phys. {123} (1989), 177-254.

[O1]
\name{E.M. Opdam},
{  Some applications of hypergeometric shift
operators}, Invent.Math.{  98} (1989), 1--18.

[O2]
\bibline, {Harmonic analysis for certain representations of
graded Hecke algebras},
Preprint Math. Inst. Univ. Leiden W93-18 (1993).

\endreferences

\bye